\documentclass[12pt]{article}
\usepackage{amssymb,amsmath,amsopn,amsfonts}
\usepackage[dvips]{color}
\usepackage{colordvi,multicol}
\usepackage{epsf}
\newfam\msbfam
\font\tenmsb=msbm10 \textfont\msbfam=\tenmsb \font\sevenmsb=msbm7
\scriptfont\msbfam=\sevenmsb \font\fivemsb=msbm5
\scriptscriptfont\msbfam=\fivemsb

\def\th#1{\vspace{1mm}\noindent{\bf #1}\quad}

\def\proof{\vspace{1mm}\noindent{\it Proof}\quad}

\numberwithin{equation}{section}

\def\bc{\begin{center}}
\def\ec{\end{center}}

\def\hang{\hangindent\parindent}
\def\textindent#1{\indent\llap{\qquad #1\ \ \enspace}\ignorespaces}
\def\ref{\par\hang\textindent}

\begin{document}

\centerline{\Large\bf Large deviation principles for the stochastic
} \centerline{\Large\bf  quasi-geostrophic
equations}\footnotetext{\footnotesize Supported in part by NSFC (No.11201234),  a project funded by the PAPD of Jiangsu Higher Education Institutions and the DFG through IRTG
1132 and CRC 701.}
\date{}
\vspace{0.5 true cm}
 \centerline{Wei Liu}
 \centerline{\small
School of Mathematics and Statistics, Jiangsu Normal University, Xuzhou 221116, China}
\centerline{\small
Department of Mathematics, University of Bielefeld, D-33615 Bielefeld, Germany}  \centerline{\small E-mail: weiliu@math.uni-bielefeld.de}
 \centerline{Michael R\"{o}ckner} \centerline{\small Department of Mathematics, University of Bielefeld, D-33615 Bielefeld, Germany}
 \centerline{\small E-mail: roeckner@mathematik.uni-bielefeld.de}
\centerline{Xiang-Chan Zhu} \centerline{\small
School of Science, Beijing Jiaotong University, Beijing 100044,
China}\centerline{\small
Department of Mathematics, University of Bielefeld, D-33615 Bielefeld, Germany}
 \centerline{\small E-mail: zhuxiangchan@126.com}

\begin{abstract}
In this paper we establish the large deviation principle for the stochastic quasi-geostrophic equation
with small multiplicative noise in the subcritical case. The proof is mainly based on the  weak convergence approach. Some analogous results are also obtained for the small time asymptotics of  the stochastic quasi-geostrophic equation.

\end{abstract}

\bigbreak

\section{Introduction}
 The main aim of this work is to establish  large deviation principles for the stochastic quasi-geostrophic equation, which is an important model in geophysical fluid dynamics.  We consider the following two dimensional (2D) stochastic quasi-geostrophic equation in the periodic domain $\mathbb{T}^2=\mathbb{R}^2/(2\pi \mathbb{Z})^2$:
$$\frac{\partial \theta(t,x)}{\partial t}=-u(t,x)\cdot \nabla \theta(t,x)-\kappa (-\triangle)^\alpha \theta(t,x)+(G(\theta)\xi)(t,x)  \eqno(1.1)$$
with  initial condition $$\theta(0,x)=\theta_0(x). \eqno(1.2)$$
Here  $0<\alpha<1, \kappa>0$ are real numbers, $\theta(t,x)$ (representing the potential temperature) is a real-valued function of $t$ and $x$,
 $\xi(t,x)$ is a Gaussian random field, white noise in time and subject to the restrictions imposed below, $u$ (representing the fluid velocity) is determined by $\theta$ via the following relation:
$$u=(u_1,u_2)=(-R_2\theta,R_1\theta)=R^\bot\theta,\eqno(1.3)$$
where $R_j$ is the $j$-th periodic Riesz transform.
The case $\alpha=\frac{1}{2}$ is called the critical case, the case $\alpha>\frac{1}{2}$ subcritical and the case $\alpha<\frac{1}{2}$ supercritical.

Equation (1.1) is used to describe models arising in meteorology and oceanography. In the deterministic case ($G=0$) such equations are important models in geophysical fluid dynamics. Indeed,
they are special cases of general quasi-geostrophic approximations for atmospheric
and oceanic fluid flows with small Rossby and Ekman numbers. These models arise under
the assumptions of fast rotation, uniform stratification and uniform potential vorticity. The case $\alpha=1/2$ exhibits similar features (singularities) as the 3D Navier-Stokes equations and can therefore serve as a
 model case  for the latter. For more details about the geophysical background, see for instance \cite{CMT94,P87}. In the deterministic case, this equation has been already intensively investigated because of both its intrinsic mathematical importance and its applications in geophysical fluid dynamics (see e.g. \cite{CV06,CW99,Ju04,Ju05,KNV07,Re95} and the references therein). For example,  the global existence of weak solutions has been obtained in \cite{Re95} and one very  remarkable result in \cite{CV06} proved the existence of a classical solution for $\alpha=\frac{1}{2}$ and the other in \cite{KNV07} proved that  solutions for $\alpha=\frac{1}{2}$ with periodic $C^\infty$ data remain $C^\infty$ for all time.

Recently, in \cite{RZZ12} the  two last named authors and Rongchan Zhu have studied the 2D stochastic quasi-geostrophic equation on
$\mathbb{T}^2$ for general parameter $\alpha\in (0,1)$ and for both additive as well as multiplicative noise. For the subcritical case $\alpha>\frac{1}{2}$ the authors obtained a (probabilistically strong) solution.
In this paper, we want to establish the large deviation principles for stochastic quasi-geostrophic equation  both for small noise and for short time in the subcritical case.

The large deviation theory  concerns the asymptotic behavior of a family of random variables $\{\theta_\varepsilon\}$ and we refer to the monographs \cite{DZ92,St84} for many historical remarks and extensive references. It asserts that for some tail or extreme event $A$, $P(\theta_\varepsilon\in A)$ converges to zero exponentially fast as $\varepsilon\rightarrow0$ and  the exact rate of convergence is given by the so-called rate function. The large deviation principle was first established by Varadhan in \cite{Va66} and  he also studied the small time asymptotics of finite dimensional diffusion processes in \cite{Va67}. Since then, many important results concerning the large deviation principle have been established. For  results on the large deviation principle for  stochastic differential equations in finite dimensional case we refer to \cite{FW84}. For the extensions to infinite dimensional diffusions or SPDE, we refer the readers to \cite{BDM08,CM10, DM08,Li09,MSS09,RZ08,SS06,SZ11,XZ09} and the references therein.

  The large deviation principle for the stochastic quasi-geostrophic equation with small multiplicative noise is proved in  Section 3 and the small time large deviations for this equation in Section 4 in the subcritical case (i.e. $\alpha>\frac{1}{2}$).
The proof of small noise LDP  is mainly based on the weak convergence approach from \cite{BD00}. Compared to  some recent works on LDP for SPDE (cf.\cite{CM10,Li09,RZ08}), the main difficulty here lies in dealing with the nonlinear term in (1.1) since the solution to the stochastic quasi-geostrophic equation is not as regular as in the case of SPDE within the variational framework (see \cite{CM10,Li09,RZ08} for many examples). For example, for 2D Navier-Stokes equation, the solution lies in the first order Sobolev space by which the nonlinear term can be dominated. Compared with this, the solution of the stochastic quasi-geostrophic equation only lies in $H^\alpha$ (see definition below) and the nonlinear term cannot be handled as for 2D Navier-Stokes equation. Here we use the regularity of solutions of the deterministic equation to control the nonlinear term. Indeed, the solution of the deterministic quasi-geostrophic equation will be in $H^\delta$ if the initial value lies in $H^\delta$ (see Theorem A.1). Our main result on small noise large deviations for equation (1.1) is formulated in Theorem 3.9.  The small time large deviation principle describes the behavior of the temperature of the fluid when  time is very small. The proof is mainly inspired by the approach from \cite{XZ09}. We first establish the large deviation principle on $L^\infty([0,T],H)$ if the initial value is smooth (see Theorem 4.1). However, since the solution to the stochastic quasi-geostrophic equation is very irregular, we cannot approximate the initial value similarly as in \cite{XZ09} for the 2D Navier-Stokes equation to obtain the result for more general initial value. In order to overcome this difficulty, we  establish the small time large deviation principle with general initial value on a larger state space (see Theorem 4.2). Here we use the $L^p$-norm estimate to control the nonlinear term. But these $L^p$-norm estimates we cannot prove  by Galerkin approximation, instead we use another approximation which can be seen as a piecewise linear equation on small subintervals (see (4.11)).

\section{Notations and preliminaries}

 In the following, we will restrict ourselves to flows which have zero average on the torus, i.e.
$$\int_{\mathbb{T}^2}\theta dx=0.$$
Thus (1.3) can be restated as
$$u=(-\frac{\partial \psi}{\partial x_2},\frac{\partial \psi}{\partial x_1})\  \textrm { and } \ (-\triangle)^{1/2}\psi=-\theta.$$
Set
$$H=\{f\in L^2(\mathbb{T}^2):\int_{\mathbb{T}^2}f dx=0\}$$
 and let $|\cdot|$ and $\langle \cdot, \cdot\rangle$ denote the usual norm and inner product in $H$ respectively. On the periodic domain $\mathbb{T}^2$, it is well known that
 $$\{\sin (kx)|k\in \mathbb{Z}^2_+\}\cup\{\cos(kx)|k\in \mathbb{Z}^2_-\}$$
  form an eigenbasis (we denote it by $\{e_k\}$)
  of $-\triangle$  and the corresponding eigenvalues are $|k|^2$. Here
  $$ \mathbb{Z}^2_+=\{(k_1,k_2)\in \mathbb{Z}^2|k_2>0\}\cup\{(k_1,0)\in \mathbb{Z}^2|k_1>0\}, \ \mathbb{Z}^2_-=\{(k_1,k_2)\in \mathbb{Z}^2|(-k_1,-k_2)\in \mathbb{Z}^2_+ \}. $$
   Now we define
   $$\|f\|_{H^s}^2=\sum_k |k|^{2s}\langle f,e_k\rangle^2$$ and let $H^s$ denote the (Sobolev) space of all $f$ such that $\|f\|_{H^s}$ is finite.

   Set $\Lambda=(-\triangle)^{1/2}$, then we have  $$\|f\|_{H^s}=|\Lambda^s f|.$$

  By the singular integral theory of Calder\'{o}n and Zygmund (cf.\cite[Chapter 3]{St70}), for any $p\in(1,\infty)$, there exists a constant $C(p)$ such that
\begin{equation}\|u\|_{L^p}\leq C(p)\|\theta\|_{L^p}.\end{equation}

For fixed $\alpha\in(0,1)$, we define the linear operator
$$ A_\alpha: D(A_\alpha)=H^{2\alpha}(\mathbb{T}^2)\subset H\rightarrow H, \  A_\alpha u=\kappa (-\triangle)^\alpha u.$$
It is well known that $A_\alpha$ is positive definite and self-adjoint with the same eigenbasis  as that of $-\bigtriangleup$ mentioned above. We denote the eigenvalues of $A_\alpha$ by $0<\lambda_1\leq\lambda_2\leq\cdots$  and renumber the above eigenbasis correspondingly as $e_1,e_2,\cdots$.

We  first recall the following  product estimate (cf.\cite[Lemma A.4]{Re95}).

\th{Lemma 2.1} Suppose that $s>0$ and $p\in (1,\infty)$. If $f,g\in C^\infty(\mathbb{T}^2)$ ,
 then
\begin{equation}\|\Lambda^s(fg)\|_{L^p}\leq C\left(\|f\|_{L^{p_1}}\|\Lambda^sg\|_{L^{p_2}}+\|g\|_{L^{p_3}}\|\Lambda^sf\|_{L^{p_4}} \right),\end{equation}
where $p_i\in (1,\infty), i=1,...,4$ satisfy that
$$\frac{1}{p}=\frac{1}{p_1}+\frac{1}{p_2}=\frac{1}{p_3}+\frac{1}{p_4}.$$

For the reader's convenience we also recall the following standard Sobolev inequality (cf.\cite[Chapter V]{St70}):

\th{Lemma 2.2} Suppose that $q>1, p\in [q,\infty)$ and
$$\frac{1}{p}+\frac{\sigma}{2}=\frac{1}{q}.$$
If $\Lambda^\sigma f\in L^q$, then we have $f\in L^p$ and there is a constant $C\geq 0$ (independent of $f$) such that
$$\|f\|_{L^p}\leq C\|\Lambda^\sigma f\|_{L^q}.$$

\section{Freidlin-Wentzell's large deviations in the subcritical case}

In this section, we consider the large deviation principle for the
stochastic quasi-geostrophic equation with small multiplicative
noise. Here we will use the weak convergence approach introduced by
Budhiraja and Dupuis in \cite{BD00}. Let us first recall some
standard definitions and results from  large deviation theory
(cf.\cite{DZ93}).

Let $\{X^\varepsilon\}$ be a family of random variables defined on a probability space $(\Omega,\mathcal{F},P)$  taking values in some Polish space $E$.

\th{Definition 3.1}(Rate function) A function $I: E\rightarrow[0,\infty]$ is called a rate function if $I$ is lower semicontinuous. A rate function $I$ is called a good rate function if the level set $\{x\in E: I(x)\leq M\}$ is compact for each $M<\infty$.

\th{Definition 3.2}(I)(Large deviation principle) The sequence $\{X^\varepsilon\}$ is said to satisfy the large deviation principle with rate function $I$ if for each Borel subset $A$ of $E$
$$-\inf_{x\in A^o}I(x)\leq\liminf_{\varepsilon\rightarrow0}\varepsilon \log P(X^\varepsilon\in A)\leq\limsup_{\varepsilon\rightarrow0}\varepsilon\log P(X^\varepsilon\in A)\leq-\inf_{x\in \bar{A}}I(x),$$
where $A^o$ and $\bar{A}$ denote  the interior and  closure of $A$ in $E$ respectively.

(II)(Laplace principle) The sequence $\{X^\varepsilon\}$ is said to satisfy the Laplace principle with rate function $I$  if for each bounded continuous real-valued function $h$ defined on $E$
$$\lim_{\varepsilon\rightarrow0}\varepsilon\log E\{\exp[-\frac{1}{\varepsilon}h(X^\varepsilon)]\}=-\inf_{x\in E}\{h(x)+I(x)\}.$$

It is well known that the large deviation principle and the Laplace principle are  equivalent
 if $E$ is a Polish space and the rate function is good. The equivalence
  is essentially a consequence of Varadhan's lemma and
Bryc's converse theorem (cf.\cite{DZ93}).

Suppose $W(t)$ is a cylindrical Wiener process on a Hilbert space $U$ (with  inner product $\langle \cdot,\cdot\rangle_U$ and norm $|\cdot|_U$) defined on a probability space $(\Omega,\mathcal{F},\mathcal{F}_t,P)$ (i.e. the paths of $W$ take values in $C([0,T],Y)$, where $Y$ is another Hilbert space such that the embedding $U\subset Y$ is Hilbert-Schmidt).
Now we define
$$\aligned
\mathcal{A}&=\left\{\phi: \phi~ \text{is a}~ U\text{-valued}~ \{\mathcal{F}_t\}\text{-predictable process s.t.} \int_0^T|\phi(s)|^2_U ds<\infty \ a.s.\right\};\\
S_M&=\left\{v\in L^2([0,T], U): \int^T_0|v(s)|^2_U ds\leq M \right\}; \\
\mathcal{A}_M&=\left\{\phi\in \mathcal{A}: \phi(\omega)\in S_M, P\text{-}a.s.\right\}.
\endaligned$$
Here we remark that we will always refer
to the weak topology on the set $S_M$ in this paper.

 Suppose $g^\varepsilon: C([0,T],Y)\rightarrow E$ is a measurable map and $X^\varepsilon=g^\varepsilon(W)$. 
 Now we formulate the following sufficient conditions for the Laplace principle (equivalently,
large deviation principle) of $X^\varepsilon$ as $\varepsilon\rightarrow0$.

\th{Hypothesis 3.3} There exists a measurable map $g^0:C([0,T],Y)\rightarrow E$ such that the following conditions hold:

1) Let $\{v^\varepsilon:\varepsilon>0\}\subset \mathcal{A}_M$ for some $M<\infty$. If $v^\varepsilon$ converge to $v$ as $S_M$-valued random elements in distribution, then $g^\varepsilon(W(\cdot)+\frac{1}{\sqrt{\varepsilon}}\int_0^\cdot v^\varepsilon(s)ds)$ converge in distribution to $g^0(\int_0^\cdot v(s)ds)$.

2) For every $M<\infty$, the set $K_M=\{g^0(\int_0^\cdot v(s)ds):v\in S_M\}$ is a compact subset of $E$.
\vskip.10in

The following crucial result was proven in \cite{BD00} (see also
\cite{BD98} for finite dimensional case).

\th{Theorem 3.4}(\cite[Theorem 4.4]{BD00})  If $\{g^\varepsilon\}$ satisfies Hypothesis 3.3, then $\{X^\varepsilon\}$ satisfies the Laplace principle (hence large deviation principle) on $E$ with the good rate function $I$ given by
$$I(f)=\inf_{\{v\in L^2([0,T],U): ~ f=g^0(\int_0^\cdot v(s)ds)\}}\left\{\frac{1}{2}\int_0^T|v(s)|_U^2ds\right\}.\eqno(3.1)$$

Now we reformulate (1.1)-(1.3) in the following form of an abstract stochastic evolution equation:
$$\left\{\begin{array}{ll}d\theta(t)+A_\alpha\theta(t)dt+u(t)\cdot \nabla\theta(t)dt= G(\theta)dW(t),&\ \ \ \ \textrm{ }\\\theta(0)=\theta_0\in H,&\ \ \ \ \textrm{ } \end{array}\right.\eqno(3.2)$$
where $u$ satisfies (1.3).
\vskip.10in

We first need to impose some assumptions on $G$ such that (3.2) has a unique solution.
Let $L_2(U,H)$ be the space of all Hilbert-Schmidt operators from $U$ to $H$ and  $\{f_n\}$ be an ONB of $U$.
Recall that we only consider the subcritical case (i.e. $\alpha> \frac{1}{2}$) in this work. Let $\beta>3$ be some fixed constant.

\th{Hypothesis 3.5}Suppose that  $G$ satisfies the following conditions:

i) There exist some positive real numbers $C_1,C_2, C_3$ and
$\rho_1<2\kappa$ such that
 $$\|G(\theta)\|^2_{L_2(U,H)}\leq C_1|\theta|^2+\rho_1|\Lambda^\alpha\theta|^2+C_2, \theta\in H^\alpha;$$
 $$\|G(\theta)\|^2_{L_2(U,H^{-\beta})}\leq C_3(|\theta|^2+1), \theta\in H^\alpha.$$

ii) If $\theta_n, \theta\in H^\alpha$ and $\theta_n\rightarrow
\theta$ in $H$, then for all $v\in C^\infty(\mathbb{T}^2)$,
$$\lim_{n\rightarrow\infty} |G(\theta_n)^*(v)-G(\theta)^*(v) |_U=0,$$
 where the asterisk denotes the adjoint operator.

iii)  For some $p$ with $0< 1/p<\alpha-\frac{1}{2}$, there exists some constant $C$ such that
 $$\int_{\mathbb{T}^2}(\sum_j|G(\theta)(f_j)|^2)^{p/2}dx\leq C\left(\int_{\mathbb{T}^2} |\theta|^pdx+1 \right),
 \  \theta\in H^\alpha\cap L^{p}(\mathbb{T}^2);
 \eqno(3.3)$$

iv) There exist some constants $C$ and $\beta_1<2\kappa$ such that
$$\| \Lambda^{-1/2}(G(\theta_1)-G(\theta_2))\|^2_{L_2(U,H)}\leq C|\Lambda^{-1/2}(\theta_1-\theta_2)|^2+\beta_1|\Lambda^{\alpha-\frac{1}{2}}(\theta_1-\theta_2)|^2,
\ \theta_1,\theta_2\in H^\alpha. \eqno(3.4)$$
\vskip.10in
Now we give the definition of the (probabilistically) strong
solution to (3.2).

\th{Definition 3.6}  We say that there exists a (probabilistically)
strong solution to (3.2) on  $[0,T]$ if for every probability space
$(\Omega,\mathcal{F},\{\mathcal{F}_t\}_{t\in [0,T]},P)$ with an
$\mathcal{F}_t$-cylindrical Wiener process $W$, there exists an
$\mathcal{F}_t$-adapted process $\theta:[0,T]\times
\Omega\rightarrow H$ such that for $P$-$a.s.$ $\omega\in \Omega$
$$\theta(\cdot,\omega)\in L^\infty([0,T];H)\cap L^2([0,T];H^\alpha)\cap C([0,T];H^{-\beta})$$
and $P$-$a.s.$
$$\langle \theta(t),\varphi\rangle+\int_0^t\langle A_\alpha^{1/2}\theta(s),A_\alpha^{1/2}\varphi\rangle ds-\int_0^t \langle u(s)
\cdot \nabla \varphi,\theta(s)\rangle ds=\langle
\theta_0,\varphi\rangle+\langle \int_0^tG(\theta(s))dW(s),
\varphi\rangle$$ for all $t\in [0,T]$ and all $\varphi\in
C^1(\mathbb{T}^2)$.
\vskip.10in \th{Remark}Note that $div u=0$, so for regular functions
$\theta$ and $\varphi$  we have
$$\langle u(s)\cdot \nabla (\theta(s)+\varphi),\theta(s)+\varphi\rangle=0.$$
Hence,
$$\langle u(s)\cdot \nabla \theta(s),\varphi \rangle=-\langle u(s)\cdot \nabla \varphi,\theta(s)\rangle.$$
This relation justifies  the integral equation in Definition 3.6.

\vskip.10in

We recall the following existence and uniqueness result from
\cite{RZZ12}.
\vskip.10in \th{Theorem 3.7}(\cite[Theorem 4.3]{RZZ12})
Assume $\alpha>\frac{1}{2}$ and Hypothesis 3.5 hold. Then for each
initial condition $\theta_0\in L^p$, there exists a pathwise unique
probabilistically strong solution $\theta$ of equation (3.2) on
$[0,T]$ with  initial condition $\theta(0)=\theta_0$ such that
$$E\sup_{t\in [0,T]}|\Lambda^{-1/2}\theta(t)|^2<\infty.$$ Moreover,
the solution $\theta$ satisfies $$ E\sup_{t\in
[0,T]}\|\theta(t)\|^p_{L^p}+
 E\int_0^T|\Lambda^\alpha\theta(t)|^2dt<\infty.$$
\vskip.10in

Now we consider the stochastic quasi-geostrophic equation driven by
small multiplicative noise:
$$d\theta^\varepsilon(t)+A_\alpha\theta^\varepsilon(t)dt+u^\varepsilon(t)\cdot \nabla\theta^\varepsilon(t)dt= \sqrt{\varepsilon}G(\theta^\varepsilon)dW(t)\eqno(3.5)$$
with $\theta^\varepsilon(0)=\theta_0\in L^p$. Here $u^\varepsilon$ satisfies (1.3) with $\theta$ replaced by $\theta^\varepsilon$. By Theorem 3.7, under
Hypothesis 3.5, there exists a pathwise unique strong solution of
(3.5) in $L^\infty([0,T],H)\cap L^2([0,T],H^\alpha)\cap
C([0,T],H^{-\beta})$. Therefore,  there exist Borel-measurable
functions
$$g^\varepsilon: C([0,T],Y)\rightarrow L^\infty([0,T],H)\cap
L^2([0,T],H^\alpha)\cap C([0,T],H^{-\beta})$$ such that
$\theta^\varepsilon(\cdot)=g^\varepsilon(W(\cdot))$. \vskip.10in Now
the aim  is to prove the large deviation principle for
$\theta^\varepsilon$. For this purpose we need to impose some
further  assumptions on $G$.

\th{Hypothesis 3.8} Assume $G$ satisfies the following conditions:

i) $G(\theta)$ is a bounded operator from $U$ to $H^\delta$ for some $\delta>2-2\alpha$ such that
 $$\|G(\theta)\|_{L(U,H^{\delta})}\leq C(\|\theta\|_{H^{\delta+\alpha}}+1), \ \theta\in H^{\delta+\alpha}\eqno(3.6)$$
 and for $r:=(2-2\alpha)\vee \alpha$  $$\|G(\theta)\|_{L(U,H^r)}\leq C(\|\theta\|_{H^{\delta+\alpha}}+1), \ \theta\in H^{\delta+\alpha}.\eqno(3.7)$$

ii) $$\|G(\theta_1)-G(\theta_2)\|_{L(U,H)}\leq C\|\theta_1-\theta_2\|_{H^\alpha},\ \theta_1,\theta_2\in H^\alpha.$$

 \th{Remark}(i) (3.6) can also be replaced by the following condition:
 $$\|G(\theta)\|_{L(U,H^{\delta-\alpha})}\leq C(\|\theta\|_{H^{\delta}}+1).$$
 (ii) Typical examples for $G$ satisfying Hypothesis 3.5 and 3.8 have the following form: for $\theta\in H^\alpha$
$$G(\theta)y=\sum_{k=1}^\infty b_k\langle y,f_k\rangle_U g(\theta), y\in U,$$
where $g\in C_b^1(\mathbb{R})$  and $b_k$ are $C^\infty$ functions on $\mathbb{T}^2$ satisfying
$$\sum_{k=1}^\infty b_k^2(\xi)\leq C, ~~~~ \sum_{k=1}^\infty |\Lambda^{\delta\vee r}b_k|^2\leq C.  $$

 \vskip.10in

For $v\in L^2([0,T],U)$, we consider the following skeleton equation
$$ \frac{d\theta_v(t)}{d t}= -A_\alpha \theta_v(t)- u_v(t)\cdot \nabla\theta_v(t)+ G(\theta_v)v(t) \eqno(3.8)$$
with $\theta_v(0)=\theta_0\in H^\delta\cap L^p$. Here $u_v$ satisfies (1.3) with $\theta$ replaced by $\theta_v$. Then by Hypothesis 3.5 and 3.8  we have
$$\|G(\theta)v\|_{L^p}\leq C|v|_U(\|\theta\|_{L^p}+1);\eqno(A.1)$$
$$ \|G(\theta)v\|_{H^{\delta}}\leq C|v|_U(\|\theta\|_{H^{\delta+\alpha}}+1);\eqno(A.2)$$
$$ | \Lambda^{-1/2}(G(\theta_1)-G(\theta_2))v|\leq |v|_U(C|\Lambda^{-1/2}(\theta_1-\theta_2)|+\sqrt{\beta_1}|\Lambda^{\alpha-\frac{1}{2}}(\theta_1-\theta_2)|).\eqno(A.3)$$
 By a similar argument as in \cite[Theorems 3.5 and 3.7]{Re95}, we know that (3.8) has a unique solution
 $\theta_v\in L^\infty([0,T],H^\delta\cap L^p)\cap L^2([0,T],H^{\delta+\alpha})\cap C([0,T],H^{-\beta})$.
 For the completeness we include the proof of this result in the Appendix.

 \th{Remark} Here we want to emphasize that although by Theorem A.1 in Appendix if $\theta_0\in H^\delta\cap L^p$, then we have $\theta_v\in L^\infty([0,T],H^\delta\cap L^p)\cap L^2([0,T],H^{\delta+\alpha})\cap C([0,T],H^{-\beta})$. However,  this might be not true for $\theta^\varepsilon$.
  This is the reason why we establish the large deviation principle for $\theta^\varepsilon$  on $L^\infty([0,T],H)\cap L^2([0,T],H^\alpha)\cap C([0,T],H^{-\beta})$ (which is the state space of $\theta^\varepsilon$) instead of $L^\infty([0,T],H^\delta)\cap L^2([0,T],H^{\delta+\alpha})\cap C([0,T],H^{-\beta})$.
 \vskip.10in
Define $g^0:C([0,T],Y)\rightarrow L^\infty([0,T],H)\cap L^2([0,T],H^\alpha)\cap C([0,T],H^{-\beta})$ by
$$g^0(h)=\left\{\begin{array}{ll}\theta_v, \textrm{ if } h=\int_0^\cdot v(s)ds \textrm{ for some } v\in L^2([0,T],U),&\ \ \ \ \textrm{ }\\0, \textrm{ otherwise. }&\ \ \ \ \textrm{ } \end{array}\right.$$

Now we formulate the main result concerning the large deviation principle for $\theta^\varepsilon$.

\th{Theorem 3.9} Suppose that Hypothesis 3.5 and Hypothesis 3.8 hold.  Then for any $\theta_0\in H^\delta\cap L^p$ with $p$ in Hypothesis 3.5 iii),
$\{\theta^\varepsilon\}$ satisfies the Laplace principle (hence
large deviation principle) on $L^\infty([0,T],H)\cap
L^2([0,T],H^\alpha)\cap C([0,T],H^{-\beta})$ with a good rate
function given by (3.1).

\proof  To prove the theorem it suffices to verify the two
conditions in Hypothesis 3.3 so that Theorem 3.4 is applicable to
obtain the large deviation principle for $\theta^\varepsilon$.
\vskip.10in

[Step 1] First we show that the set $K_M=\{g^0(\int_0^\cdot v(s)ds):v\in S_M\}$ is a compact subset of $L^\infty([0,T],H)\cap L^2([0,T],H^\alpha)\cap C([0,T],H^{-\beta})$.

Let $\{\theta_n\}$ be a sequence in $K_M$ where $\theta_n$ corresponds to the solution of (3.8) with $v_n\in S_M$ in place of $v$. By the weak compactness of $S_M$ in $L^2([0,T],U)$, there exists a subsequence (which we still denote it by $\{v_n\}$)  converging to a limit $v$ weakly in $L^2([0,T],U)$.

Let $w_n=\theta_n-\theta_v$, it suffices to show that
$w_n\rightarrow 0$ (in fact, a subsequence is enough) in
$L^\infty([0,T],H)\cap L^2([0,T],H^\alpha)\cap C([0,T],H^{-\beta})$
as $n\rightarrow \infty$.

Note that $u_n\cdot\nabla w_n\in H^{-\alpha}$, where $u_n$ satisfies (1.3) with $\theta$ replaced by $\theta_n$. In fact, we have uniform $L^p$ norm bound for $\theta_n,w_n$ by
Theorem A.1. And we also have
$$\aligned |{}_{H^{-\alpha}}\!\langle u_n\cdot\nabla w_n, \varphi\rangle_{H^\alpha}|&=|{}_{H^{-\alpha}}\!\langle \nabla\cdot (u_nw_n), \varphi\rangle_{H^\alpha}|\leq |\Lambda^\alpha\varphi||\Lambda^{1-\alpha}(u_n\cdot w_n)|\\&\leq |\Lambda^\alpha\varphi|(|\Lambda^{1-\alpha+\sigma}w_n|\|\theta_n\|_{L^p}+|\Lambda^{1-\alpha+\sigma}\theta_n|\|w_n\|_{L^p}),\endaligned$$
where $\sigma=\frac{2}{p}<2\alpha-1$ and  we use $div u_n=0$ in the first equality and Lemmas 2.1, 2.2 and (2.1) in the last inequality.
Thus by $div u_n=0$, we obtain
$${}_{H^{-\alpha}}\!\langle u_n\cdot\nabla w_n, w_n\rangle_{H^\alpha} =0.\eqno(3.9)$$

If $\delta<1$ we get
$$\aligned |\langle (u_n-u_v)\cdot\nabla\theta_v, w_n\rangle |&=|\langle \nabla\cdot((u_n-u_v)\theta_v), w_n\rangle |\leq |\Lambda^\alpha w_n||\Lambda^{1-\alpha}((u_n-u_v)\cdot\theta_v)|\\&\leq C|\Lambda^\alpha w_n|(|\Lambda^{2-\alpha-\delta}w_n||\Lambda^{\delta}\theta_v|+|\Lambda^{1-\alpha+\delta-(1-\alpha)}\theta_v||\Lambda^{2-\alpha-\delta         }w_n|)\\&\leq
C|\Lambda^\alpha w_n||\Lambda^{\alpha}w_n|^\gamma
|w_n|^{1-\gamma}|\Lambda^{\delta}\theta_v|\\&\leq
\frac{\kappa}{4}|\Lambda^\alpha
w_n|^2+C|\Lambda^\delta\theta_v|^{N}|w_n|^2,
\endaligned\eqno(3.10)$$
where $\gamma=\frac{2-\alpha-\delta}{\alpha},
N=\frac{2\alpha}{2\alpha-2+\delta}$ and we use $div (u_n-u_v)=0$ in the first equality, Lemmas 2.1, 2.2 and (2.1) in the second inequality, the interpolation
inequality and $\delta>2-2\alpha$ in the third inequality and Young's inequality in the last
inequality.

Similarly, if $\delta\geq1$ we get
$$\aligned |\langle (u_n-u_v)\cdot\nabla\theta_v, w_n\rangle |&=|\langle \nabla\cdot((u_n-u_v)\theta_v), w_n\rangle |\leq |\Lambda^\alpha w_n||\Lambda^{1-\alpha}((u_n-u_v)\cdot\theta_v)|\\&\leq C|\Lambda^\alpha w_n||\Lambda^{1-\alpha+\sigma_1}w_n||\Lambda^{\delta}\theta_v|\\&\leq
C|\Lambda^\alpha
w_n||\Lambda^{\alpha}w_n|^{\gamma_1}|w_n|^{1-\gamma_1}|\Lambda^{\delta}\theta_v|\\&\leq
\frac{\kappa}{4}|\Lambda^\alpha
w_n|^2+C|\Lambda^\delta\theta_v|^{N_1}|w_n|^2,
\endaligned$$
where $0<\sigma_1<2\alpha-1,
\gamma_1=\frac{1-\alpha+\sigma_1}{\alpha},N_1=\frac{2\alpha}{2\alpha-1-\sigma_1}$ and we use $div (u_n-u_v)=0$ in the first equality, Lemmas 2.1, 2.2 and (2.1),  $\delta\geq1$ in the second inequality, the interpolation
inequality in the third inequality and Young's inequality in the last
inequality.

In the following we only prove the result for $\delta<1$ and the
argument for $\delta\geq1$ is similar.

By (3.8) we have
 $$\aligned
 ~ & |w_n(t)|^2+2\kappa\int_0^t |\Lambda^\alpha w_n|^2 ds\\
 =& 2\int_0^t \left(-\langle u_n\cdot\nabla\theta_n, w_n\rangle+\langle u_v\cdot\nabla\theta_v,
 w_n\rangle \right) ds\\&+2\int_0^t\langle G(\theta_n)v_n-G(\theta_v)v,w_n\rangle ds\\
 = & -2\int_0^t\langle (u_n-u_v)\cdot\nabla\theta_v, w_n\rangle ds\\&+2\int_0^t\langle (G(\theta_n)-G(\theta_v))v_n ,w_n \rangle ds\\
 &+2\int_0^t\langle G(\theta_v)(v_n -v ),w_n \rangle ds\\
 \leq & \int_0^t\bigg[ \kappa |\Lambda^\alpha w_n|^2+C(|\Lambda^\delta\theta_v|^{N}+|v_n|_U^2)|w_n|^2\\&+2\langle G(\theta_v)(v_n -v ),w_n \rangle \bigg] ds,\endaligned$$
 where in the second equality we use (3.9) and in the last inequality we use (3.10), Hypothesis 3.8  ii) and Young's inequality.

 Let $$h_n(t):=\int_0^tG(\theta_v)(v_n-v)ds,$$
 then we have
 $$\aligned\sup_{t\in[0,T]}\|P_kh_n(t)-h_n(t)\|_{H^r}&\leq \int_0^T\|(P_k-I)G(\theta_v)\|_{L(U,H^r)}\|v_n-v\|_Udt\\&\leq (2M)^{1/2}\left(\int_0^T\|(P_k-I)G(\theta_v)\|_{L(U,H^r)}^2dt \right)^{1/2}\rightarrow0\ \textrm{ as } k\rightarrow\infty.\endaligned$$
 Here $P_k$ is the orthogonal projection in $H$ onto the space spanned by $e_1,...e_k$ and we use (3.7) and $\theta_v\in L^2([0,T]; H^{\delta+\alpha})$ which follows from Theorem A.1 in the last step.

 Since $P_kH^r\subset H^r$ is compact and $v_n\rightarrow v$ weakly in $L^2([0,T];U)$, by (3.7) it
 is easy to show that $P_kh_n\rightarrow 0$ in $C([0,T],H^r)$ as $n\rightarrow\infty$ (see e.g. \cite[Lemma 3.2]{Li09})
 using the Arz\`{e}la-Ascoli theorem  (since for any subsequence the limit is the same, this convergence holds for the whole sequence).
 Hence we obtain that $h_n\rightarrow0$ in $C([0,T],H^r)$ as $n\rightarrow\infty$.

 And we also have
 $$\aligned  ~ & \int_0^t\langle G(\theta_v)(v_n(s)-v(s)),w_n(s)\rangle ds\\
 =&\langle w_n(t),h_n(t)\rangle-\int_0^t\langle w_n'(s),h_n(s)\rangle ds\\
 =&
 \langle w_n(t),h_n(t)\rangle+\int_0^t\langle A_\alpha w_n+u_n\cdot\nabla\theta_n-u_v\cdot\nabla\theta_v,h_n\rangle ds\\&-\int_0^t\langle G(\theta_n)v_n-G(\theta_v)v,h_n\rangle ds\\=:&I_1+I_2+I_3.\endaligned\eqno(3.11)$$
 Note that
$$I_1\leq \varepsilon|w_n(t)|^2+C|h_n(t)|^2;$$
and by Hypothesis 3.5 i) and (A.4)
$$\aligned I_3&\leq \sup_{s\in [0,T]}|h_n(s)|\int_0^T \left(\|G(\theta_n)\|_{L_2(U,H)}|v_n|_U+\|G(\theta_v)\|_{L_2(U,H)}|v|_U \right)ds\\ &\leq C\sup_{s\in [0,T]}\|h_n(s)\|_{H^r}\big(\int_0^T(|\Lambda^\alpha\theta_v|^2+|\Lambda^\alpha\theta_n|^2+C)ds\big)^{1/2}\leq C\sup_{s\in [0,T]}\|h_n(s)\|_{H^r}.\endaligned$$
For $\varphi\in H^{2-2\alpha}$, we obtain
 $$\aligned |\langle u_n\cdot\nabla\theta_n-u_v\cdot\nabla\theta_v,\varphi\rangle|&=|\langle \nabla \cdot (u_n\theta_n-u_v\theta_v),\varphi\rangle|
 \\&\leq C|\Lambda^{2\alpha-1}(u_n\theta_n-u_v\theta_v)||\Lambda^{2-2\alpha}\varphi|
 \\&\leq C(|\Lambda^\alpha \theta_n|^2+|\Lambda^\alpha \theta_v|^2)|\Lambda^{2-2\alpha}\varphi|,\endaligned$$
 where we use $div u_n=0$ and $div u_v=0$ in the first equality and Lemmas 2.1, 2.2 and (2.1) in the last inequality.

Hence $$\|u_n\cdot\nabla\theta_n-u_v\cdot\nabla\theta_v\|_{H^{-(2-2\alpha)}}\leq
 C\left( |\Lambda^\alpha \theta_n|^2+|\Lambda^\alpha \theta_v|^2 \right) .$$
 Therefore,
$$\aligned I_2\leq & \int_0^t(\|A_\alpha w_n(s)\|_{H^{-\alpha}}+\|u_n\cdot\nabla\theta_n-u_v\cdot\nabla\theta_v\|_{H^{-(2-2\alpha)}})\|h_n(s)\|_{H^r}ds\\\leq &C\sup_{s\in [0,T]}\|h_n(s)\|_{H^r}\int_0^t(\|w_n\|_{H^\alpha}+\|\theta_n\|^2_{H^\alpha}+\|\theta_v\|^2_{H^\alpha})ds\\
\leq &C\sup_{s\in [0,T]}\|h_n(s)\|_{H^r},\endaligned$$where in the last step we use (A.4).

Then the Gronwall lemma and (3.11) yield that
$$\sup_{t\in[0,T]}|w_n(t)|^2+\frac{\kappa}{2}\int_0^T |\Lambda^\alpha w_n|^2 ds\leq C\sup_{t\in [0,T]}\|h_n(t)\|_{H^r} \left(\exp\left\{C\int_0^T \left(|\Lambda^\delta \theta_v|^{N}+|v_n|_U^2 \right) ds \right\}+1\right).$$
 Then by (A.4) we have
$$\sup_{t\in [0,T]}|w_n(t)|^2+\frac{\kappa}{2}\int_0^T |\Lambda^\alpha w_n|^2 ds\rightarrow 0,\ \ n\rightarrow\infty.$$

[Step 2]  Suppose that $\{v_\varepsilon: \varepsilon>0 \} \subset \mathcal{A}_M  $ for some $M<\infty$ and   $v_\varepsilon$ converge to $v$  as $S_M$-valued random elements in distribution.
 Then, by Girsanov's theorem,  $\bar{\theta}_{v_\varepsilon}=g^\varepsilon(W(\cdot)+\frac{1}{\sqrt{\varepsilon}}\int_0^\cdot v^\varepsilon(s)ds)$ solves the following equation
$$d\bar{\theta}_{v_\varepsilon}(t)+A_\alpha \bar{\theta}_{v_\varepsilon}(t)dt+u_{\bar{\theta}_{v_\varepsilon}}(t)\cdot \nabla\bar{\theta}_{v_\varepsilon}(t)dt= G(\bar{\theta}_{v_\varepsilon})v_\varepsilon(t)dt+\sqrt{\varepsilon}G(\bar{\theta}_{v_\varepsilon})dW(t).\eqno(3.12)$$
Here $u_{\bar{\theta}_{v_\varepsilon}}$ satisfies (1.3) with
$\theta$ replaced by $\bar{\theta}_{v_\varepsilon}$.

Since $S_M$ is a Polish space, by the Skorohod theorem, we can construct processes $(\tilde{v}_\varepsilon,\tilde{v},\tilde{W}_\varepsilon)$ such that the joint distribution of $(\tilde{v}_\varepsilon,\tilde{W}_\varepsilon)$ is the same as that of $(v_\varepsilon,W)$,  the distribution of $v$ coincides with that of $\tilde{v}$  and $\tilde{v}_\varepsilon\rightarrow \tilde{v}$ a.s. as $S_M$-valued random elements.

Setting  $w_\varepsilon:=\bar{\theta}_{\tilde{v}_\varepsilon}-\theta_{\tilde{v}}$, it suffices to prove that $w_\varepsilon\rightarrow 0$ in probability in $L^\infty([0,T],H)\cap L^2([0,T],H^\alpha)\cap C([0,T],H^{-\beta})$. For $w_\varepsilon$ and $u_{\bar{\theta}_{v_\varepsilon}}$ we also have similar estimates as (3.9) and (3.10). In the following we write $v_\varepsilon=\tilde{v}_\varepsilon, W=\tilde{W}_\varepsilon$ for simplicity.

It\^{o}'s formula and (3.9) imply that
$$\aligned
~ & |w_\varepsilon(t)|^2+2\kappa\int_0^t |\Lambda^\alpha w_\varepsilon|^2 ds\\
=& 2\int_0^t \left( -\langle u_{\bar{\theta}_{v_\varepsilon}}\cdot\nabla\bar{\theta}_{v_\varepsilon}, w_\varepsilon\rangle+\langle u_v\cdot\nabla\theta_v, w_\varepsilon\rangle \right) ds\\
&+2\int_0^t\langle G(\bar{\theta}_{v_\varepsilon}(s))v_\varepsilon(s)-G(\theta_v(s))v(s),w_\varepsilon(s)\rangle ds\\
&+2\sqrt{\varepsilon}\int_0^t\langle w_\varepsilon,G(\bar{\theta}_{v_\varepsilon})dW\rangle
+\varepsilon\int_0^t\|G(\bar{\theta}_{v_\varepsilon})\|^2_{L_2(U,H)}ds
\\= & -2\int_0^t\langle (u_{\bar{\theta}_{v_\varepsilon}}-u_v)\cdot\nabla\theta_v, w_\varepsilon\rangle ds +2\int_0^t\langle ( G(\bar{\theta}_{v_\varepsilon}(s))-G(\theta_v(s)))v_\varepsilon(s),w_\varepsilon(s)\rangle ds\\&+2\int_0^t\langle G(\theta_v(s))(v_\varepsilon(s)-v(s)),w_\varepsilon(s)\rangle ds\\&+2\sqrt{\varepsilon}\int_0^t\langle w_\varepsilon,G(\bar{\theta}_{v_\varepsilon})dW\rangle
+\varepsilon\int_0^t\|G(\bar{\theta}_{v_\varepsilon})\|^2_{L_2(U,H)}ds\\
\leq & \int_0^t \left[ \kappa |\Lambda^\alpha w_\varepsilon|^2+C(|\Lambda^\delta\theta_v|^{N}+|v_\varepsilon|^2_U)|w_\varepsilon|^2 \right] ds\\
&+2\int_0^t\langle G(\theta_v(s))(v_\varepsilon(s)-v(s)),w_\varepsilon(s)\rangle ds\\&+2\sqrt{\varepsilon}\int_0^t\langle w_\varepsilon,G(\bar{\theta}_{v_\varepsilon})dW\rangle
+\varepsilon\int_0^t\|G(\bar{\theta}_{v_\varepsilon})\|^2_{L_2(U,H)}ds,\endaligned\eqno(3.13)$$
where in the last inequality we use (3.10), Hypothesis 3.8 ii) and Young's inequality.

Similarly we define
$$h_\varepsilon(t)=\int_0^t G(\theta_v(s))(v_\varepsilon(s)-v(s))ds.$$
 Then by the same argument as [Step 1] we know  $h_\varepsilon(t)\rightarrow 0$ in $C([0,T],H^r)$ a.s. as $\varepsilon\rightarrow0$.

 By It\^{o}'s formula and a similar argument as in (3.11) we have
$$\aligned
~ &\int_0^t\langle G(\theta_v(s))(v_\varepsilon(s)-v(s)),w_\varepsilon(s)\rangle ds\\
\leq & \varepsilon |w_\varepsilon(t)|^2+C\left(1+\int_0^t|\Lambda^\alpha\bar{\theta}_{v_\varepsilon}|^2ds \right) \sup_{s\in [0,T]}\|h_\varepsilon(s)\|_{H^r}-\sqrt{\varepsilon}\int_0^t\langle h_\varepsilon,G(\bar{\theta}_{v_\varepsilon})dW\rangle.\endaligned$$
Define $$\tau_{L,\varepsilon}:=T\wedge\inf\{t:|\bar{\theta}_{v_\varepsilon}(t)|^2+\int_0^t|\Lambda^\alpha\bar{\theta}_{v_\varepsilon}(s)|^2ds>L\}.$$
Since $\bar{\theta}_{v_\varepsilon}$ is weakly continuous in $H$, $\tau_{L,\varepsilon}$ is a stopping time with respect to $\mathcal{F}_{t+}=\cap_{s>t}\mathcal{F}_s$ and $|\bar{\theta}_{v_\varepsilon}(t\wedge\tau_{L,\varepsilon})|\leq L$.
By the Burkholder-Davis-Gundy inequality one has
$$\aligned
~ &\sqrt{\varepsilon}E\sup_{t\in [0,\tau_{L,\varepsilon}]}|\int_0^t\langle w_\varepsilon-h_\varepsilon,G(\theta_{v_\varepsilon})dW\rangle| \\
\leq & C\sqrt{\varepsilon}E\left(\int_0^{\tau_{L,\varepsilon}}|\bar{\theta}_{v_\varepsilon}-\theta_v-h_\varepsilon|^2
\|G(\bar{\theta}_{v_\varepsilon})\|^2_{L_2(U,H)}ds \right)^{1/2}\\\leq & C\sqrt{\varepsilon}E\left(\int_0^{\tau_{L,\varepsilon}}
(|\Lambda^\alpha\bar{\theta}_{v_\varepsilon}|^2+1)ds \right)^{1/2}
\leq C\sqrt{\varepsilon}.
\endaligned$$
Combining the above estimates with (3.13) and applying Gronwall's lemma we have
$$\aligned
~ &\sup_{s\in [0,t]}|w_\varepsilon(s)|^2+\frac{\kappa}{2}\int_0^t|\Lambda^\alpha w_\varepsilon|^2 ds\\
\leq &\bigg[C(1+\int_0^t|\Lambda^\alpha\bar{\theta}_{v_\varepsilon}|^2ds)\sup_{s\in [0,T]}\|h_\varepsilon(s)\|_{H^r} +2\sqrt{\varepsilon}\sup_{t\in [0,T]}|\int_0^t\langle w_\varepsilon-h_\varepsilon,G(\bar{\theta}_{v_\varepsilon})dW\rangle|\\
&+\varepsilon\int_0^t\|G(\bar{\theta}_{v_\varepsilon})\|^2_{L_2(U,H)}ds\bigg]
\exp\left\{C\int_0^T \left(|\Lambda^\delta \theta_v|^{N}+|v_\varepsilon|_U^2 \right)dr\right\}.\endaligned$$

Then we have $$\sup_{t\in [0,\tau_{L,\varepsilon}]}|w_\varepsilon(t)|^2+\frac{\kappa}{2}\int_0^{\tau_{L,\varepsilon}}|\Lambda^\alpha w_\varepsilon|^2ds\rightarrow 0$$ in probability as $\varepsilon\rightarrow0$.

 By It\^{o}'s formula and standard argument (cf. \cite[Theorem 3.3]{RZZ12}) we have
$$\sup_{\varepsilon\in[0,1)}E[\sup_{t\in[0,T]}|\bar{\theta}_{v_\varepsilon}(t)|^2+\int_0^T|\Lambda^\alpha\bar{\theta}_{v_\varepsilon}(t)|^2dt]<\infty.$$
Let $L$ be fixed. Then for a suitable constant $C$
$$\sup_{\varepsilon\in[0,1)} P(\tau_{L,\varepsilon}=T)\geq 1-\frac{C}{L}.$$
Therefore,  we have
$$\sup_{t\in [0,T]}|w_\varepsilon(t)|^2+\frac{\kappa}{2}\int_0^T |\Lambda^\alpha w_\varepsilon|^2ds\rightarrow 0$$ in probability as $\varepsilon\rightarrow0$.

Now the proof of Theorem 3.9 is complete.
$\hfill\Box$
\section{The small time large deviations  in the subcritical case}

In this section, we establish some  small time large deviations
results for the stochastic quasi-geostrophic equation.  The proof  is
mainly inspired by  the approach used in \cite{XZ09}. We  consider
the stochastic quasi-geostrophic equation (3.2) again and assume
that $G$ satisfies Hypothesis 3.5. Then by Theorem 3.7, for
$\theta_0\in L^p$ there exists a pathwise unique strong solution of
(3.2) in $L^\infty([0,T],H)\cap L^2([0,T],H^\alpha)\cap
C([0,T],H^{-\beta})$ for  $\beta>3$.

\vskip.10in

We assume the following additional  conditions on $G$:

S.1) There exists a constant $L$ such that for some $\delta>0$
$$\|G(\theta)\|_{L_2(U,H^\delta)}^2\leq L(1+\|\theta\|^2_{H^\delta}), \  \theta\in H^\delta;$$

S.2) There exists a constant $L_1$ such that
$$\|G(\theta_1)-G(\theta_2)\|_{L_2(U,H)}^2\leq L_1|\theta_1-\theta_2|^2, \
\theta_1,\theta_2\in H. $$
\vskip.10in

Let $\varepsilon>0$. By the scaling property of the Wiener process, it is easy to see that $\theta(\varepsilon t)$ coincides in law with the solution of the following equation
$$d\theta^\varepsilon(t)+\varepsilon A_\alpha\theta^\varepsilon(t)dt+\varepsilon u^\varepsilon(t)\cdot \nabla\theta^\varepsilon(t)dt= \sqrt{\varepsilon}G(\theta^\varepsilon)dW(t)\eqno(4.1)$$
with $\theta^\varepsilon(0)=\theta_0$. Here $u^\varepsilon$ satisfies (1.3) with $\theta$ replaced by $\theta^\varepsilon$.

Let $\mu^\varepsilon$ be the law of $\theta^\varepsilon$ on $L^\infty([0,T], H)$. Now we formulate  the  small time large deviation principle for (4.1) on $L^\infty([0,T], H)$ for regular initial value $\theta_0$.

\th{Theorem 4.1} Suppose that S.1) for some $\delta\geq\alpha$ and $\delta>2-2\alpha$, S.2) and
Hypothesis 3.5 hold.
Then for $\theta_0\in H^\delta\cap L^p$ with $p$ in Hypothesis 3.5 iii), $\mu^\varepsilon$ satisfies the
large deviation principle on $L^\infty([0,T], H)$ with  rate
function $I$ given by $$I(f)=\inf_{\{v\in
L^2([0,T],U):f(t)=\theta_0+\int_0^tG(f(s))v(s)ds\}}
\left\{\frac{1}{2}\int_0^T|v(s)|_U^2ds \right\}.\eqno(4.2)$$

\proof
Let $v^\varepsilon$ be the solution of the stochastic equation
$$v^\varepsilon(t)=\theta_0+\sqrt{\varepsilon}\int_0^tG(v^\varepsilon(s))dW(s) \eqno(4.3)$$
 and $\nu^\varepsilon$ be the law of $v^\varepsilon$ on $L^\infty([0,T], H)$. Then by \cite{Li09} we know that $\nu^\varepsilon$ satisfies the large  deviation principle with  rate function $I$ given by (4.2).
Now it is sufficient to show that the two families of  probability measures $\mu^\varepsilon$ and $\nu^\varepsilon$ are exponentially equivalent, i.e. for any $\eta>0$,
$$\lim_{\varepsilon\rightarrow0}\varepsilon\log P(\sup_{0\leq t\leq T}| \theta^\varepsilon(t)-v^\varepsilon(t) |^2>\eta)=-\infty.\eqno(4.4)$$
Then the conclusion in Theorem 4.1 follows directly from \cite[Theorem 4.2.13]{DZ93}.

In the following we assume that $\delta<1$ and for $\delta\geq1$ the proof is similar.

For $M>0$, we define the following stopping times:
$$\tau_{\varepsilon,M}=\inf\{t\ge 0: \|v^\varepsilon(t)\|_{H^\delta}^2>M\}.$$
 Then we have
$$\aligned &P(\sup_{0\leq t\leq T}|\theta^\varepsilon(t)-v^\varepsilon(t)|^2>\eta,\sup_{0\leq t\leq T}\|v^\varepsilon(t)\|^2_{H^\delta}\leq M)\\\leq & P(\sup_{0\leq t\leq T\wedge \tau_{\varepsilon,M} }|\theta^\varepsilon(t)-v^\varepsilon(t)|^2>\eta).\endaligned\eqno(4.5)$$
Applying It\^{o}'s formula to $|v^\varepsilon(t\wedge \tau_{\varepsilon,M})-\theta^\varepsilon(t\wedge \tau_{\varepsilon,M} )|^2$ we get
$$\aligned & |v^\varepsilon(t\wedge \tau_{\varepsilon,M})-\theta^\varepsilon(t\wedge \tau_{\varepsilon,M} )|^2+
2\varepsilon\kappa\int_0^{t\wedge\tau_{\varepsilon,M}} |\Lambda^{\alpha}(v^\varepsilon(s)-\theta^\varepsilon(s))|^2ds
\\=&2\varepsilon\int_0^{t\wedge\tau^n_{\varepsilon,M}}\langle A_\alpha v^\varepsilon(s),(v^\varepsilon(s)-\theta^\varepsilon(s))\rangle ds
+2\varepsilon\int_0^{t\wedge\tau_{\varepsilon,M}}\langle u^\varepsilon\cdot\nabla\theta^\varepsilon,(v^\varepsilon-\theta^\varepsilon)\rangle ds
\\&+2\sqrt{\varepsilon}\int_0^{t\wedge\tau_{\varepsilon,M}}\langle v^\varepsilon-\theta^\varepsilon,(G(v^\varepsilon)-G(\theta^\varepsilon))dW\rangle
\\&+\varepsilon\int_0^{t\wedge\tau_{\varepsilon,M}}\|G(v^\varepsilon)-G(\theta^\varepsilon)\|^2_{L_2(U,H)}ds.\endaligned$$
Note that by  similar arguments as in (3.9) and (3.10), we have
$$\aligned  & |\langle u^\varepsilon\cdot\nabla\theta^\varepsilon,v^\varepsilon-\theta^\varepsilon\rangle|\\
=&|\langle u^\varepsilon\cdot\nabla(\theta^\varepsilon-v^\varepsilon),
\theta^\varepsilon-v^\varepsilon\rangle\\&+\langle (u^\varepsilon-u_{v^\varepsilon})\cdot\nabla v^\varepsilon,
\theta^\varepsilon-v^\varepsilon\rangle\\&+\langle u_{v^\varepsilon}\cdot\nabla v^\varepsilon,\theta^\varepsilon-v^\varepsilon\rangle|\\\leq& \frac{\kappa}{2}|\Lambda^{\alpha}(\theta^\varepsilon-v^\varepsilon)|^2+ C|\Lambda^\delta v^\varepsilon|^N|\theta^\varepsilon-v^\varepsilon|^2+C|\Lambda^{\delta}v^\varepsilon|^4,
\endaligned$$
where $u_{v_\varepsilon}$ satisfies (1.3) with $\theta$ replaced by $v_\varepsilon$.
Here for the last term we use the following estimate:
$$\aligned |\langle u_{v^\varepsilon}\cdot\nabla v^\varepsilon,\theta^\varepsilon-v^\varepsilon\rangle|&=|\langle\nabla\cdot (u_{v^\varepsilon}v^\varepsilon), \theta^\varepsilon-v^\varepsilon\rangle|\\
&\leq C|\Lambda^\alpha(\theta^\varepsilon-v^\varepsilon)||\Lambda^{1-\alpha}(u_{v_\varepsilon}v^\varepsilon)|\\&\leq|\Lambda^{\alpha}(\theta_n^\varepsilon-v_n^\varepsilon)|
|\Lambda^{\delta}v^\varepsilon|^2,
\endaligned$$
where in the first equality we use $div u_{v_\varepsilon}=0$ and in the last inequality we use Lemmas 2.1, 2.2 and $\delta\geq (2-2\alpha)\vee\alpha$.

Therefore, by S.2) and Young's inequality we get
$$\aligned &|v^\varepsilon(t\wedge \tau_{\varepsilon,M})-\theta^\varepsilon(t\wedge \tau_{\varepsilon,M} )|^2+
2\varepsilon\kappa\int_0^{t\wedge\tau_{\varepsilon,M}} |\Lambda^{\alpha}(v^\varepsilon(s)-\theta^\varepsilon(s))|^2ds
\\\leq &2\varepsilon\int_0^{t\wedge\tau_{\varepsilon,M}} \left( \frac{\kappa}{2}|\Lambda^{\alpha}(v^\varepsilon-\theta^\varepsilon)|^2
+C|\Lambda^{\alpha}v^\varepsilon|^2 \right)ds
\\&+2\varepsilon\int_0^{t\wedge\tau_{\varepsilon,M}} \left(\frac{\kappa}{2}|\Lambda^{\alpha}(\theta^\varepsilon-v^\varepsilon)|^2+ C|\Lambda^\delta v^\varepsilon|^N|\theta^\varepsilon-v^\varepsilon|^2+C|\Lambda^\delta v^\varepsilon|^4 \right)ds
\\&+2\sqrt{\varepsilon}\int_0^{t\wedge\tau_{\varepsilon,M}}\langle v^\varepsilon-\theta^\varepsilon,(G(v^\varepsilon)-G(\theta^\varepsilon))dW\rangle
\\&+\varepsilon C\int_0^{t\wedge\tau_{\varepsilon,M}}|v^\varepsilon-\theta^\varepsilon|^2ds.\endaligned$$
Then by Gronwall's lemma and $\delta\geq \alpha$ we have
$$\aligned &|v^\varepsilon(t\wedge \tau_{\varepsilon,M})-\theta^\varepsilon(t\wedge \tau_{\varepsilon,M} )|^2
\leq \bigg[2\varepsilon\int_0^{t\wedge\tau_{\varepsilon,M}} (C|\Lambda^{\delta}v^\varepsilon|^2+C|\Lambda^\delta v^\varepsilon|^4) ds\\
 &+ 2\sqrt{\varepsilon}|\int_0^{t\wedge\tau_{\varepsilon,M}}\langle v^\varepsilon-\theta^\varepsilon,(G(v^\varepsilon)
 -G(\theta^\varepsilon))dW\rangle|\bigg]e^{\varepsilon C\int_0^{t\wedge\tau_{\varepsilon,M}}|\Lambda^\delta v^\varepsilon|^Nds+Ct\varepsilon}
.\endaligned$$
To estimate the stochastic integral term, we will use the following result from \cite{BY82,Da76}, namely that there exists a universal constant $c$ such that for any $q\geq2$ and for any continuous martingale $M_t$ with $M_0=0$, one has
$$\|M^*_t\|_{L^q}\leq cq^{1/2}\|\langle M\rangle_t^{1/2}\|_{L^q},\eqno(4.6)$$
where $M^*_t=\sup_{0\leq s\leq t} |M_s|$.

By this result and S.2) we have
$$\aligned &(E[\sup_{0\leq s\leq t\wedge\tau_{\varepsilon,M}} |v^\varepsilon(s)-\theta^\varepsilon(s)|^2]^q)^{2/q}\\\leq&Ce^{\varepsilon CM^{N/2}t+Ct\varepsilon}\left[(\varepsilon Mt+\varepsilon M^2t)^2+q\varepsilon \big(E(\int_0^{t\wedge\tau_{\varepsilon,M}} |v^\varepsilon(r)-\theta^\varepsilon(r)|^4dr)^{q/2}\big)^{2/q} \right]\\\leq&Ce^{\varepsilon CM^{N/2}t+Ct\varepsilon}\left[(\varepsilon Mt+\varepsilon M^2t)^2+q\varepsilon\int_0^t(E[\sup_{0\leq r\leq s\wedge\tau_{\varepsilon,M}} |v^\varepsilon(r)-\theta^\varepsilon(r)|^2]^q)^{2/q}ds \right].\endaligned$$
Then Gronwall's lemma yields that
$$\aligned &(E[\sup_{0\leq s\leq T\wedge\tau_{\varepsilon,M}} |v^\varepsilon(s)-\theta^\varepsilon(s)|^2]^q)^{2/q}\\\leq&Ce^{\varepsilon CM^{N/2}T+CT\varepsilon}(\varepsilon MT+\varepsilon M^2T)^2\exp\left[ CqT\varepsilon e^{\varepsilon CM^{N/2}T+CT\varepsilon} \right].\endaligned$$
Fixing $M$ and taking $q=2/\varepsilon$ we obtain
$$\aligned &\varepsilon\log P(\sup_{0\leq t\leq T\wedge\tau_{\varepsilon,M}}|\theta^\varepsilon(t)-v^\varepsilon(t)|^2>\eta)\\
\leq &\varepsilon\log \frac{E[\sup_{0\leq s\leq T\wedge\tau_{\varepsilon,M}} |v^\varepsilon(s)-\theta^\varepsilon(s)|^{2q}]}{\eta^q}\\\leq & \log C(\varepsilon MT+\varepsilon M^2T)^2-2\log\eta+CTe^{\varepsilon CM^{N/2}T+CT\varepsilon}+\varepsilon CM^{N/2}T+CT\varepsilon\\
\rightarrow& -\infty,\textrm{ as }\varepsilon\rightarrow0.\endaligned$$
Therefore, by (4.5) there exists  $\varepsilon_0$ such that for every $\varepsilon$ satisfying $0<\varepsilon\leq\varepsilon_0$,
$$P(\sup_{0\leq t\leq T}|\theta^\varepsilon(t)-v^\varepsilon(t)|^2>\eta,
\sup_{0\leq t\leq T}\|v^\varepsilon(t)\|^2_{H^\delta}\leq M)\leq e^{-R/\varepsilon}.\eqno(4.7)$$
By the same argument as in \cite[Lemma 3.2]{XZ09} and  S.1) we have
 $$\lim_{M\rightarrow\infty}\sup_{0<\varepsilon\leq1}\varepsilon\log P(\sup_{0\leq t\leq T}\|v^\varepsilon(t)\|^2_{H^\delta}>M)=-\infty.\eqno(4.8)$$
 Then for any $R>0$, there exists a constant $M$ such that for every $\varepsilon\in(0,1]$ the following inequality holds:
$$P(\sup_{0\leq t\leq T}\| v^\varepsilon(t)\|_{H^\delta}^2>M)\leq e^{-R/\varepsilon}.\eqno(4.9)$$
 By (4.7) and (4.9), we know that there exists  $\varepsilon_0$ such that for every $\varepsilon$ satisfying $0<\varepsilon\leq\varepsilon_0$ we have
$$P(\sup_{0\leq t\leq T}|\theta^\varepsilon(t)-v^\varepsilon(t)|^2>\eta)\leq2e^{-R/\varepsilon}.$$
Since $R$ is arbitrary, we obtain (4.4).

Hence the proof of Theorem 4.1 is complete.
$\hfill\Box$

\vskip.10in

 Note that the solution of (4.1) is not as regular as in the case of the 2D stochastic Navier-Stokes equation. In Theorem 4.1 we use the regularity of $v^\varepsilon$ to control the nonlinear term, but we can not approximate the initial value in (4.1) to obtain the large deviation principle on $L^\infty([0,T],H)$ for general initial value in $L^p$ as  Xu and Zhang did in \cite{XZ09} for the 2D stochastic Navier-Stokes equation since the nonlinear term can not be dominated.
To overcome this difficulty, now we enlarge the state space of the solution and  use $L^p$ norm estimate to control the nonlinear term. Then we establish the large deviation principle on $L^\infty([0,T], H^{-1/2})$.
\vskip.10in
We consider the following condition on $G$.

S.3) There exists a constant $L_2$ such that
$$\|\Lambda^{-1/2}(G(\theta_1)-G(\theta_2))\|_{L_2(U,H)}^2\leq L_2|\Lambda^{-1/2}(\theta_1-\theta_2)|^2, \
\theta_1,\theta_2\in H^\alpha. $$

\th{Remark} Typical examples for $G$ satisfying Hypothesis 3.5 and S.1)-S.3) have the following form: for $\theta\in H^\alpha$
$$G(\theta)y=\sum_{k=1}^\infty b_k\langle y,f_k\rangle_U \theta,  y\in U,$$
 where $b_k$ are $C^\infty$ functions on $\mathbb{T}^2$ satisfying $\sum_{k=1}^\infty b_k^2(\xi)\leq M$ and $\sum_{k=1}^\infty |\Lambda^{1+\varepsilon}b_k|^2\leq M$ for some $\varepsilon>0$.

\vskip.10in
Let $\bar{\mu}^\varepsilon$ be the law of $\theta^\varepsilon$ on $L^\infty([0,T], H^{-1/2})$. Now we formulate our main result about the  small time large deviation principle for (4.1).

\th{Theorem 4.2} Suppose that S.1) for $\delta\geq(\frac{3}{4}-\frac{\alpha}{2})\vee (\alpha-\frac{1}{2})$ and Hypothesis 3.5, S.3) hold.  Then for $\theta_0\in L^p$, $\bar{\mu}^\varepsilon$ satisfies the large deviation principle on $L^\infty([0,T], H^{-1/2})$ with  rate function $I$ given by (4.2).

 It is sufficient to show that the two families of  probability measures $\bar{\mu}^\varepsilon$ and $\nu^\varepsilon$ (for simplicity we still use the same notation) are exponentially equivalent, i.e. for any $\eta>0$,
$$\lim_{\varepsilon\rightarrow0}\varepsilon\log P(\sup_{0\leq t\leq T}|\Lambda^{-1/2}(\theta^\varepsilon(t)-v^\varepsilon(t))|^2>\eta)=-\infty.\eqno(4.10)$$
Then the conclusion in Theorem 4.2 follows directly from \cite[Theorem 4.2.13]{DZ93}.

In order to show (4.10)  we prove a few lemmas in below.

\th{Lemma 4.3} $$\lim_{M\rightarrow\infty}\sup_{0<\varepsilon\leq1}\varepsilon\log P(\sup_{0\leq t\leq T}\|\theta^\varepsilon(t)\|^p_{L^p}>M)=-\infty.$$

\proof We consider the same approximation $\theta^{\varepsilon,n}$ to $\theta^\varepsilon$ as in \cite[Theorem 3.3]{RZZ12}. We pick a smooth function $\phi\geq0$ such that supp $\phi\subset[1,2]$ and  $\int_0^\infty \phi=1$. Then for $\sigma>0$ we define
$$U_\sigma[\theta](t):=\int_0^\infty \phi(\tau)(k_\sigma*R^\bot\theta)(t-\sigma\tau)d\tau,$$ where $k_\sigma$ is the periodic Poisson kernel in $\mathbb{T}^2$ given by $\widehat{k_\sigma}(\zeta)=e^{-\sigma|\zeta|},\zeta\in \mathbb{Z}^2$, and we set $\theta(t)=0$ for $t<0$.

We take a sequence $\delta_n$ converging to $0$ and consider the following  equation:
$$d\theta^{\varepsilon,n}(t)+\varepsilon A_\alpha \theta^{\varepsilon,n}(t)dt+\varepsilon u^{\varepsilon,n}(t)\cdot \nabla\theta^{\varepsilon,n}(t)dt= \sqrt{\varepsilon}k_{\delta_n}*G(\theta^{\varepsilon,n})dW(t) \eqno(4.11)$$ with initial data $\theta^{\varepsilon,n}(0)=k_{\delta_n}*\theta_0$ and $u^{\varepsilon,n}=U_{\delta_n}[\theta^{\varepsilon,n}]$. For a fixed $n$, this is a linear equation in $\theta^{n,\varepsilon}$ on each subinterval $[t^n_k,t^n_{k+1}]$ with $t^n_k=k\delta_n$, since $u^{\varepsilon,n}$ is determined by the values of $\theta^{\varepsilon,n}$ on the two previous subintervals.

Then by \cite[Theorem 3.3, Step 2]{RZZ12} , there exists a weak solution to (4.11) which converges in distribution to $\theta^\varepsilon$ in $L^2([0,T],H)\cap C([0,T],H^{-\beta})$.

 By \cite[Lemma 5.1]{Kr10} we have (here we write  $\theta(t)=\theta^{\varepsilon,n}(t), u(t)=u^{\varepsilon,n}(t)$
 to simplify the notation)
$$\aligned
\|\theta(t)\|_{L^p}^p=&\|k_{\delta_n}*\theta_0\|_{L^p}^p+\varepsilon\int_0^t\bigg[-p\int_{\mathbb{T}^2}|\theta(s)|^{p-2}\theta(s) (\Lambda^{2\alpha}\theta(s)+u(s)\cdot\nabla\theta(s))dx\\&
+\frac{1}{2}p(p-1)\varepsilon\int_{\mathbb{T}^2}|\theta(s)|^{p-2}(\sum_j|k_{\delta_n}*G(\theta(s))(f_j)|^2)dx\bigg]ds
\\ &+p\sqrt{\varepsilon}\int_0^t
\int_{\mathbb{T}^2}|\theta(s)|^{p-2}\theta(s) k_{\delta_n}*G(\theta(s))dxdW(s)\\
\leq& \|\theta_0\|_{L^p}^p+\frac{1}{2}p(p-1)\varepsilon\int_0^t
\int_{\mathbb{T}^2}|\theta(s)|^{p-2}(\sum_j|k_{\delta_n}*G(\theta(s))(f_j)|^2)dxds
\\&+p\sqrt{\varepsilon}\int_0^t
\int_{\mathbb{T}^2}|\theta(s)|^{p-2}\theta(s)k_{\delta_n}*G(\theta(s))dx dW(s)\\
\leq&\|\theta_0\|_{L^p}^p+\varepsilon\int_0^t\left(\int_{\mathbb{T}^2}|\theta(s)|^pdx
+C\int(\sum_j|k_{\delta_n}*G(\theta(s))(f_j)|^2)^{p/2}dx \right)ds
\\&+p\sqrt{\varepsilon}\int_0^t\int_{\mathbb{T}^2}|\theta(s)|^{p-2}\theta(s) k_{\delta_n}*G(\theta(s))dxdW(s),\endaligned$$
where in the first inequality we used $div u=0$ and $\int|\theta|^{p-2}\theta\Lambda^{2\alpha}\theta\geq0$ (cf. \cite[Lemma 3.2]{Re95}) as well as Young's inequality in the second inequality.

Then by Hypothesis 3.5 (iii) we have$$\aligned\sup_{t\in[0,T]}\|\theta(t)\|_{L^p}^p\leq&\|\theta_0\|_{L^p}^p+\varepsilon CT+C\varepsilon\int_0^T\sup_{t\in[0,s]}\|\theta(t)\|_{L^p}^pds
\\+&p\sqrt{\varepsilon}\sup_{0\leq t\leq T}|\int_0^t\int_{\mathbb{T}^2}|\theta(s)|^{p-2}\theta(s) k_{\delta_n}*G(\theta(s))dxdW(s)|.\endaligned$$
Therefore, for $q\geq2$ we obtain
$$\aligned (E(\sup_{t\in[0,T]}\|\theta(t)\|_{L^p}^{pq}))^{1/q}\leq&\|\theta_0\|_{L^p}^p+\varepsilon CT+C\varepsilon (E(\int_0^T\sup_{t\in[0,s]}\|\theta(t)\|_{L^p}^pds)^q)^{1/q}
\\+&p\sqrt{\varepsilon}(E\sup_{0\leq t\leq T}|\int_0^t\int_{\mathbb{T}^2}|\theta(s)|^{p-2}\theta(s) k_{\delta_n}*G(\theta(s))dxdW(s)|^q)^{1/q}.\endaligned$$

Using (4.6) and Minkowski's inequality we have
$$\aligned
&~~  p\sqrt{\varepsilon}(E\sup_{0\leq t\leq T}|\int_0^t\int_{\mathbb{T}^2}|\theta(s)|^{p-2}\theta(s) k_{\delta_n}*G(\theta(s))dxdW(s)|^q)^{1/q}\\
&\leq pc\sqrt{q\varepsilon}(E(\int_0^T(\int_{\mathbb{T}^2}|\theta(s)|^{p-1}(\sum_j|k_{\delta_n}*G(\theta(s))(f_j)|^2)^{1/2}dx)^2ds)^{q/2})^{1/q}
\\ &\leq pc\sqrt{q\varepsilon}(E(\sup_{s\in[0,T]}\|\theta(s)\|_{L^p}^{p-1}(\int_0^T(\int_{\mathbb{T}^2}(\sum_j|k_{\delta_n}*G(\theta(s))(f_j)|^2)^{p/2}dx)^{2/p}ds)^{1/2})^q)^{1/q}
 \\&\leq pc\sqrt{q\varepsilon}(E(\sup_{s\in[0,T]}\|\theta(s)\|_{L^p}^{p-1}(\int_0^T(\int_{\mathbb{T}^2}(\sum_j|k_{\delta_n}*G(\theta(s))(f_j)|^2)^{p/2}dx)ds)^{1/p})^q)^{1/q}\\&\leq \frac{1}{2}(E\sup_{s\in[0,T]}\|\theta(s)\|_{L^p}^{pq})^{1/q}+c(p)(q\varepsilon)^{p/2}
 (E(\int_0^T(\int_{\mathbb{T}^2}(\sum_j|k_{\delta_n}*G(\theta(s))(f_j)|^2)^{p/2}dx)ds)^{q})^{1/q}\\
 &\leq \frac{1}{2}(E\sup_{s\in[0,T]}\|\theta(s)\|_{L^p}^{pq})^{1/q}+c(p)(q\varepsilon)^{p/2}
 \left[\int_0^T \left(1+(E\|\theta(s)\|^{pq}_{L^p})^{1/q} \right)ds \right], \endaligned$$
where in the last inequality we use Hypothesis 3.5 iii) and Jesen's inequality.

Hence,
$$\aligned (E(\sup_{t\in[0,T]}\|\theta(t)\|_{L^p}^{pq}))^{1/q}\leq&2\|\theta_0\|_{L^p}^p+\varepsilon CT+C\varepsilon \int_0^T(E\sup_{t\in[0,s]}\|\theta(t)\|_{L^p}^{pq})^{1/q}ds\\+&c(p)(q\varepsilon)^{p/2}
\left[\int_0^T \left(1+(E\|\theta(s)\|^{pq}_{L^p})^{1/q} \right)ds \right].\endaligned$$
Applying Gronwall's lemma  we obtain that
$$(E(\sup_{t\in[0,T]}\|\theta(t)\|_{L^p}^{pq}))^{1/q}\leq \left[2\|\theta_0\|_{L^p}^p+\varepsilon CT+c(p)(q\varepsilon)^{p/2}T \right]\exp\left[CT\varepsilon+c(p)T(q\varepsilon)^{p/2} \right].$$
Letting $n\rightarrow\infty$ we get $$(E(\sup_{t\in[0,T]}\|\theta^\varepsilon(t)\|_{L^p}^{pq}))^{1/q}\leq\left[2\|\theta_0\|_{L^p}^p+\varepsilon CT+c(p)(q\varepsilon)^{p/2}T \right]\exp\left[CT\varepsilon+c(p)T(q\varepsilon)^{p/2} \right].$$
Since $$P(\sup_{0\leq t\leq T}\|\theta^\varepsilon(t)\|^p_{L^p}>M)\leq M^{-q}E(\sup_{t\in[0,T]}\|\theta^\varepsilon(t)\|_{L^p}^{pq}),$$ letting $q=2/\varepsilon$ we get
$$\aligned &\varepsilon\log P(\sup_{0\leq t\leq T}\|\theta^\varepsilon(t)\|^p_{L^p}>M)\leq-2\log M+2\log(E(\sup_{t\in[0,T]}\|\theta^\varepsilon(t)\|_{L^p}^{pq}))^{1/q}\\\leq &-2\log M+2\log(2\|\theta_0\|_{L^p}^p+\varepsilon CT+CT)+2CT\varepsilon+2CT.\endaligned$$
Hence the proof is complete.$\hfill\Box$
\vskip.10in
Since $H^\delta\cap L^p$ is dense in $L^p$, there exists a sequence $\{\theta_0^n\} \subset H^\delta\cap L^p$ such that $$\lim_{n\rightarrow\infty} \|\theta_0^n-\theta_0\|_{L^p}=0. $$
Let $\theta_n^\varepsilon$ be the solution of (4.1) with  initial value $\theta_0^n$. From the proof of Lemma 4.3, it follows that

$$\lim_{M\rightarrow\infty}\sup_n\sup_{0<\varepsilon\leq1}\varepsilon\log P(\sup_{0\leq t\leq T}\|\theta_n^\varepsilon(t)\|^p_{L^p}>M)=-\infty.\eqno(4.12)$$

Let $v_n^\varepsilon$ be the solution of (4.3) with  initial value $\theta_0^n$. By the same argument as in (4.8) and Lemma 4.3 we have the following result.
 \vskip.10in
\th{Lemma 4.4} For every $n\in \mathbb{Z}^+$, $$\lim_{M\rightarrow\infty}\sup_{0<\varepsilon\leq1}\varepsilon\log P(\sup_{0\leq t\leq T}(\|v_n^\varepsilon(t)\|^2_{H^\delta}+\|v_n^\varepsilon(t)\|_{L^p}^p)>M)=-\infty.$$
\vskip.10in
\th{Lemma 4.5} For every $\eta>0$, $$\lim_{n\rightarrow\infty}\sup_{0<\varepsilon\leq1}\varepsilon\log P(\sup_{0\leq t\leq T}\|\theta_n^\varepsilon(t)-\theta^\varepsilon(t)\|^2_{H^{-1/2}}>\eta)=-\infty.$$

\proof For $M>0$, we define the following stopping time for $N_0=\frac{\alpha}{\alpha-\frac{1}{2}-\frac{1}{p}}$:
$$\bar{\tau}_{\varepsilon,M}=\inf\{t\ge 0:\int_0^t\|\theta^\varepsilon(t)\|_{L^p}^{N_0}dt>M\}.$$

Clearly, we have
$$\aligned &P(\sup_{0\leq t\leq T}\|\theta_n^\varepsilon(t)-\theta^\varepsilon(t)\|^2_{H^{-1/2}}>\eta,\int_0^T\|\theta^\varepsilon(t)\|^{N_0}_{L^p}dt\leq M)\\\leq& P(\sup_{0\leq t\leq T\wedge \bar{\tau}_{\varepsilon,M} }\|\theta_n^\varepsilon(t)-\theta^\varepsilon(t)\|^2_{H^{-1/2}}>\eta).\endaligned\eqno(4.13)$$
Let $k$ be a positive constant and $N_0=\frac{\alpha}{\alpha-\frac{1}{2}-\frac{1}{p}}$.  Then
applying Ito's formula to
 $$e^{-k\varepsilon\int_0^{t\wedge \bar{\tau}_{\varepsilon,M} }\|\theta^\varepsilon(s)\|_{L^p}^{N_0}ds}|\Lambda^{-1/2}(\theta^\varepsilon(t\wedge \bar{\tau}_{\varepsilon,M})-\theta_n^\varepsilon(t\wedge \bar{\tau}_{\varepsilon,M} ))|^2 $$
 we get
$$\aligned
~&~ e^{-k\varepsilon\int_0^{t\wedge \bar{\tau}_{\varepsilon,M} }\|\theta^\varepsilon(s)\|_{L^p}^{N_0}ds}|\Lambda^{-1/2}(\theta^\varepsilon(t\wedge \bar{\tau}_{\varepsilon,M})-\theta_n^\varepsilon(t\wedge \bar{\tau}_{\varepsilon,M} ))|^2\endaligned$$
$$\aligned &+
2\varepsilon\kappa\int_0^{t\wedge\bar{\tau}_{\varepsilon,M}} e^{-k\varepsilon\int_ 0^s\|\theta^\varepsilon(r)\|_{L^p}^{N_0}dr}|\Lambda^{\alpha-\frac{1}{2}}(\theta^\varepsilon(s)
-\theta^\varepsilon_n(s))|^2ds
\\=&|\Lambda^{-\frac{1}{2}}(\theta_0-\theta_0^n)|^2
-k\varepsilon\int_0^{t\wedge\bar{\tau}_{\varepsilon,M}}e^{-k\varepsilon
\int_0^s\|\theta^\varepsilon\|^{N_0}_{L^p}dr}\|\theta^\varepsilon(s)\|^{N_0}_{L^p}
|\Lambda^{-\frac{1}{2}}(\theta^\varepsilon(s)-\theta^\varepsilon_n(s))|^2ds
\\&-2\varepsilon\int_0^{t\wedge\bar{\tau}_{\varepsilon,M}}
e^{-k\varepsilon\int_0^s\|\theta^\varepsilon\|^{N_0}_{L^p}dr}\langle u^\varepsilon(s)\cdot\nabla\theta^\varepsilon(s)-u^\varepsilon_n(s)\cdot\nabla\theta^\varepsilon_n(s),
\Lambda^{-1}(\theta^\varepsilon(s)-\theta_n^\varepsilon(s))\rangle ds
\\&+2\sqrt{\varepsilon}\int_0^{t\wedge\bar{\tau}_{\varepsilon,M}}
e^{-k\varepsilon\int_0^s\|\theta^\varepsilon\|^{N_0}_{L^p}dr}\langle \Lambda^{-1/2}(\theta^\varepsilon(s)-\theta^\varepsilon_n(s)),
\Lambda^{-1/2}(G(\theta^\varepsilon(s))-G(\theta^\varepsilon_n(s)))dW(s)\rangle
\\&+\varepsilon\int_0^{t\wedge\bar{\tau}_{\varepsilon,M}}
e^{-k\varepsilon\int_0^s\|\theta^\varepsilon\|^{N_0}_{L^p}dr}\|\Lambda^{-1/2}(G(\theta^\varepsilon(s))
-G(\theta^\varepsilon_n(s)))\|^2_{L_2(U,H)}ds,\endaligned$$
where $u_n^\varepsilon$ satisfies (1.3) with $\theta$ replaced by $\theta_n^\varepsilon$.

Note that
$$ \aligned
~& \langle u^\varepsilon\cdot\nabla\theta^\varepsilon-u^\varepsilon_n\cdot\nabla\theta^\varepsilon_n,
\Lambda^{-1}(\theta^\varepsilon-\theta_n^\varepsilon)\rangle \\
=&\langle (u_n^\varepsilon-u^\varepsilon)\cdot\nabla\theta_n^\varepsilon,\Lambda^{-1}(\theta_n^\varepsilon-\theta^\varepsilon)\rangle+\langle u^\varepsilon\cdot\nabla(\theta_n^\varepsilon-\theta^\varepsilon),
\Lambda^{-1}(\theta_n^\varepsilon-\theta^\varepsilon)\rangle.
\endaligned$$
Moreover, we also have (cf.e.g. \cite{Re95})
$$\langle (u_n^\varepsilon-u^\varepsilon)\cdot\nabla\theta_n^\varepsilon,
\Lambda^{-1}(\theta_n^\varepsilon-\theta^\varepsilon)\rangle=0\eqno(4.14)$$
and
$$\aligned & | \langle u^\varepsilon\cdot\nabla(\theta_n^\varepsilon-\theta^\varepsilon), \Lambda^{-1}(\theta_n^\varepsilon-\theta^\varepsilon)\rangle|=| \langle u^\varepsilon\cdot\nabla\Lambda^{-1}(\theta_n^\varepsilon-\theta^\varepsilon),\theta_n^\varepsilon-\theta^\varepsilon\rangle|\\
\leq & \|u^\varepsilon\|_{L^p}\|\theta_n^\varepsilon-\theta^\varepsilon\|_{L^{p'}}\|\nabla \Lambda^{-1}(\theta_n^\varepsilon-\theta^\varepsilon)\|_{L^{p'}}\\\leq& C\|u^\varepsilon\|_{L^p}\|\theta_n^\varepsilon-\theta^\varepsilon\|_{H^{1/{p}}}\|\nabla \Lambda^{-1}(\theta_n^\varepsilon-\theta^\varepsilon)\|_{H^{1/{p}}}\\\leq &C\|u^\varepsilon\|_{L^p}\|\Lambda^{-1}(\theta_n^\varepsilon-\theta^\varepsilon)\|_{H^{1+\frac{1}{p}}}^2\\\leq & C\|\theta^\varepsilon\|_{L^p}\|\Lambda^{-1}(\theta_n^\varepsilon-\theta^\varepsilon)\|^{2/N}_{H^{1/2}}\|\Lambda^{-1}(\theta_n^\varepsilon-\theta^\varepsilon)\|^{2(1-\frac{1}{N})}_{H^{\frac{1}{2}+\alpha}}\\\leq& \kappa|\Lambda^{\alpha-\frac{1}{2}}(\theta_n^\varepsilon-\theta^\varepsilon)|^2+ C_0\|\theta^\varepsilon\|_{L^p}^{N_0}|\Lambda^{-1/2}(\theta_n^\varepsilon-\theta^\varepsilon)|^2,
\endaligned\eqno(4.15)$$
where $\frac{1}{p}+\frac{2}{p'}=1$ and we  used that $div u^\varepsilon=0$ in the first equality and $H^{1/p}\hookrightarrow L^{p'}$
in the second inequality, the interpolation inequality in the forth inequality and Young's inequality in the last inequality.

Therefore, by (4.14), (4.15) and S.3)
$$\aligned
~ &e^{-k\varepsilon\int_0^{t\wedge \bar{\tau}_{\varepsilon,M} }\|\theta^\varepsilon(s)\|_{L^p}^{N_0}ds}|\Lambda^{-1/2}(\theta^\varepsilon(t\wedge \bar{\tau}_{\varepsilon,M})-\theta_n^\varepsilon(t\wedge \bar{\tau}_{\varepsilon,M} ))|^2\\
&+
2\varepsilon\kappa\int_0^{t\wedge\bar{\tau}_{\varepsilon,M}} e^{-k\varepsilon\int_0^s\|\theta^\varepsilon(r)\|_{L^p}^{N_0}dr}
|\Lambda^{\alpha-\frac{1}{2}}(\theta^\varepsilon(s)-\theta^\varepsilon_n(s))|^2ds
 \\\leq&|\Lambda^{-\frac{1}{2}}
 (\theta_0-\theta_0^n)|^2-k\varepsilon\int_0^{t\wedge\bar{\tau}_{\varepsilon,M}}
 e^{-k\varepsilon\int_0^s\|\theta^\varepsilon\|^{N_0}_{L^p}dr}\|\theta^\varepsilon(s)\|^{N_0}_{L^p}
 |\Lambda^{-\frac{1}{2}}(\theta^\varepsilon(s)-\theta^\varepsilon_n(s))|^2ds\\
 &+2\varepsilon\int_0^{t\wedge\bar{\tau}_{\varepsilon,M}}e^{-k\varepsilon\int_0^s
\|\theta^\varepsilon\|^{N_0}_{L^p}dr}(\kappa|\Lambda^{\alpha-\frac{1}{2}}
(\theta_n^\varepsilon(s)-\theta^\varepsilon(s))|^2+ C_0\|\theta^\varepsilon(s)\|_{L^p}^{N_0}|\Lambda^{-1/2}(\theta_n^\varepsilon(s)-\theta^\varepsilon(s))|^2 )ds
\\&+2\sqrt{\varepsilon}\int_0^{t\wedge\bar{\tau}_{\varepsilon,M}}
e^{-k\varepsilon\int_0^s\|\theta^\varepsilon(s)\|^{N_0}_{L^p}dr}\langle \Lambda^{-1/2}(\theta^\varepsilon(s)-\theta^\varepsilon_n(s)),
\Lambda^{-1/2}(G(\theta^\varepsilon(s))-G(\theta^\varepsilon_n(s)))dW(s)\rangle
\\&+C\varepsilon\int_0^{t\wedge\bar{\tau}_{\varepsilon,M}}e^{-k\varepsilon\int_0^s
\|\theta^\varepsilon\|^{N_0}_{L^p}dr}|\Lambda^{-1/2}(\theta^\varepsilon(s)-\theta^\varepsilon_n(s))|^2ds.\endaligned$$
Choosing $k>2C_0$ and using (4.6) we have
$$\aligned &(E[\sup_{0\leq s\leq t\wedge\bar{\tau}_{\varepsilon,M}} e^{-k\varepsilon\int_0^s \|\theta^\varepsilon(r)\|_{L^p}^{N_0}dr}|\Lambda^{-1/2}(\theta^\varepsilon(s)-\theta_n^\varepsilon(s))|^2]^q)^{2/q}\\\leq&2|\Lambda^{-\frac{1}{2}}(\theta_0-\theta_0^n)|^4\\+&C(q\varepsilon+t\varepsilon^2)\int_0^t(E[\sup_{0\leq r\leq s\wedge\bar{\tau}_{\varepsilon,M}} e^{-k\varepsilon\int_0^s \|\theta^\varepsilon(r)\|_{L^p}^{N_0}dr}|\Lambda^{-1/2}(\theta^\varepsilon(s)-\theta_n^\varepsilon(s))|^2]^q)^{2/q}ds.\endaligned$$
Applying Gronwall's lemma we obtain
$$\aligned &(E[\sup_{0\leq s\leq T\wedge\bar{\tau}_{\varepsilon,M}} e^{-k\varepsilon\int_0^s \|\theta^\varepsilon(r)\|_{L^p}^{N_0}dr}|\Lambda^{-1/2}
(\theta^\varepsilon(s)-\theta_n^\varepsilon(s))|^2]^q)^{2/q} \leq 2|\Lambda^{-\frac{1}{2}}(\theta_0-\theta_0^n)|^4e^ {CT(q\varepsilon+\varepsilon^2T)}.\endaligned$$
Hence we have
$$\aligned (E[\sup_{0\leq s\leq T\wedge\bar{\tau}_{\varepsilon,M}} |\Lambda^{-1/2}(\theta^\varepsilon(s)-\theta_n^\varepsilon(s))|^2]^q)^{2/q}
\leq 2e^{2kM}|\Lambda^{-\frac{1}{2}}(\theta_0-\theta_0^n)|^4e^ {CT(q\varepsilon+\varepsilon^2T)}.\endaligned$$
Fixing $M$ and taking $q=2/\varepsilon$ we get
$$\aligned &\sup_{0<\varepsilon\leq1}\varepsilon\log P(\sup_{0\leq t\leq T\wedge\bar{\tau}_{\varepsilon,M}}\|\theta_n^\varepsilon(t)-\theta^\varepsilon(t)\|^2_{H^{-1/2}}>\eta)\\\leq &\sup_{0<\varepsilon\leq1}\varepsilon\log \frac{E[\sup_{0\leq s\leq T\wedge\bar{\tau}_{\varepsilon,M}} |\Lambda^{-1/2}(\theta^\varepsilon(s)-\theta_n^\varepsilon(s))|^{2q}]}{\eta^q}\\\leq & 2kM+\log2|\Lambda^{-\frac{1}{2}}(\theta_0-\theta_0^n)|^4-2\log\eta+C\\
\rightarrow & -\infty,\textrm{ as }n\rightarrow\infty.\endaligned\eqno(4.16)$$
By Lemma 4.3, for any $R>0$ there exists a constant $M$ such that for every $\varepsilon\in(0,1]$ the following inequality holds:
$$P(\int_0^T\|\theta^\varepsilon(t)\|^{N_0}_{L^p}dt>M)\leq P(\sup_{0\leq t\leq T}\|\theta^\varepsilon(t)\|^p_{L^p}>(\frac{M}{T})^{p/N_0})\leq e^{-R/\varepsilon}.\eqno(4.17)$$
For such  $M$, according to (4.13) and (4.16), there exists a constant $N_2$ such that for every $n\geq N_2$,
$$\sup_{0<\varepsilon\leq1}\varepsilon \log P(\sup_{0\leq t\leq T}\|\theta_n^\varepsilon(t)-\theta^\varepsilon(t)\|_{H^{-1/2}}>\eta,\int_0^T\|\theta^\varepsilon(t)\|^{N_0}_{L^p}dt\leq M)\leq -R.\eqno(4.18)$$
Combining (4.17) and (4.18) we conclude that there exists a positive integer $N_2$ such that for every $n\geq N_2$ and $\varepsilon\in(0,1]$
$$P(\sup_{0\leq t\leq T}\|\theta_n^\varepsilon(t)-\theta^\varepsilon(t)\|^2_{H^{-1/2}}>\eta)\leq2e^{-R/\varepsilon}.$$
Since $R$ is arbitrary, the assertion of the lemma follows. $\hfill\Box$
\vskip.10in
The next lemma can be proved similarly as Lemma 4.5.

\th{Lemma 4.6} For every $\eta>0$, $$\lim_{n\rightarrow\infty}\sup_{0<\varepsilon\leq1}\varepsilon\log P(\sup_{0\leq t\leq T}\|v_n^\varepsilon(t)-v^\varepsilon(t)\|^2_{H^{-1/2}}>\eta)=-\infty.$$
\vskip.10in
\th{Lemma 4.7} For every $\eta>0$ and every positive integer $n$,
$$\lim_{\varepsilon\rightarrow0}\varepsilon\log P(\sup_{0\leq t\leq T}\|\theta_n^\varepsilon(t)-v_n^\varepsilon(t))\|_{H^{-1/2}}^2>\eta)=-\infty.$$

\proof For $M>0$, we define the following stopping times:
$$\tau^n_{\varepsilon,M}=\inf\{t:\|v_n^\varepsilon(t)\|_{H^\delta}^2+\|v_n^\varepsilon(t)\|_{L^p}^p>M\}.$$
 Then we have
$$\aligned &P(\sup_{0\leq t\leq T}\|\theta_n^\varepsilon(t)-v_n^\varepsilon(t)\|^2_{H^{-1/2}}>\eta,\sup_{0\leq t\leq T}(\|v_n^\varepsilon(t)\|^2_{H^\delta}+\|v_n^\varepsilon(t)\|_{L^p}^p)\leq M)\\\leq & P(\sup_{0\leq t\leq T\wedge \tau^n_{\varepsilon,M} }\|\theta_n^\varepsilon(t)-v_n^\varepsilon(t)\|^2_{H^{-1/2}}>\eta).\endaligned\eqno(4.19)$$
Applying It\^{o}'s formula to $|\Lambda^{-1/2}(v_n^\varepsilon(t\wedge \tau^n_{\varepsilon,M})-\theta_n^\varepsilon(t\wedge \tau^n_{\varepsilon,M} ))|^2$ we get
$$\aligned & |\Lambda^{-1/2}(v_n^\varepsilon(t\wedge \tau^n_{\varepsilon,M})-\theta_n^\varepsilon(t\wedge \tau^n_{\varepsilon,M} ))|^2+
2\varepsilon\kappa\int_0^{t\wedge\tau^n_{\varepsilon,M}} |\Lambda^{\alpha-\frac{1}{2}}(v_n^\varepsilon(s)-\theta^\varepsilon_n(s))|^2ds
\\=&2\varepsilon\int_0^{t\wedge\tau^n_{\varepsilon,M}}\langle A_\alpha v_n^\varepsilon(s),\Lambda^{-1}(v_n^\varepsilon(s)-\theta^\varepsilon_n(s))\rangle ds
+2\varepsilon\int_0^{t\wedge\tau^n_{\varepsilon,M}}\langle u^\varepsilon_n\cdot\nabla\theta^\varepsilon_n,\Lambda^{-1}(v_n^\varepsilon-\theta_n^\varepsilon)\rangle ds
\\&+2\sqrt{\varepsilon}\int_0^{t\wedge\tau^n_{\varepsilon,M}}\langle \Lambda^{-1/2}(v_n^\varepsilon-\theta^\varepsilon_n),\Lambda^{-1/2}(G(v_n^\varepsilon)-G(\theta^\varepsilon_n))dW\rangle
\\&+\varepsilon\int_0^{t\wedge\tau^n_{\varepsilon,M}}\|\Lambda^{-1/2}(G(v_n^\varepsilon)-G(\theta^\varepsilon_n))\|^2_{L_2(U,H)}ds.\endaligned$$

Note that by  similar arguments as in (4.14) and (4.15), we have
$$\aligned  & |\langle u^\varepsilon_n\cdot\nabla\theta^\varepsilon_n,\Lambda^{-1}(v_n^\varepsilon-\theta_n^\varepsilon)\rangle|\\
=&|\langle (u_n^\varepsilon-u_{v_n}^\varepsilon)\cdot\nabla\theta_n^\varepsilon,
\Lambda^{-1}(\theta_n^\varepsilon-v_n^\varepsilon)\rangle\\&+\langle u_{v_n}^\varepsilon\cdot\nabla(\theta_n^\varepsilon-v_n^\varepsilon),
\Lambda^{-1}(\theta_n^\varepsilon-v_n^\varepsilon)\rangle\\&+\langle u_{v_n}^\varepsilon\cdot\nabla v_n^\varepsilon,\Lambda^{-1}(\theta_n^\varepsilon-v_n^\varepsilon)\rangle|\\\leq& \frac{\kappa}{2}|\Lambda^{\alpha-\frac{1}{2}}(\theta_n^\varepsilon-v_n^\varepsilon)|^2+ C\|v_n^\varepsilon\|_{L^p}^{N_0}|\Lambda^{-1/2}(\theta_n^\varepsilon-v_n^\varepsilon)|^2+C\|v^\varepsilon_n\|^4_{H^\delta},\endaligned$$
where $u_{v_n}^\varepsilon$ satisfies (1.3) with $\theta$ replaced by $v_n^\varepsilon$.
Here in the last step for the last term we use the following estimate:
$$\aligned |\langle u_{v_n}^\varepsilon\cdot\nabla v_n^\varepsilon,\Lambda^{-1}(\theta_n^\varepsilon-v_n^\varepsilon)\rangle|&=|\langle u_{v_n}^\varepsilon\cdot\nabla \Lambda^{-1}(\theta_n^\varepsilon-v_n^\varepsilon),v_n^\varepsilon\rangle|\leq \|\theta_n^\varepsilon-v_n^\varepsilon\|_{L^{p_1}}\|v_n^\varepsilon\|^2_{L^{p_2}}\\&\leq|\Lambda^{\alpha-\frac{1}{2}}(\theta_n^\varepsilon-v_n^\varepsilon)|
|\Lambda^{\delta}v_n^\varepsilon|^2,
\endaligned$$
where $\frac{1}{p_1}+\frac{2}{p_2}=1,\frac{1}{p_1}+\frac{\alpha-1/2}{2}=\frac{1}{2}$ and we use $H^{\alpha-\frac{1}{2}}\subset L^{p_1}$ and $H^\delta\subset L^{p_2}$ since $\delta\geq(\frac{3}{4}-\frac{\alpha}{2})$.

Therefore, by S.3)
$$\aligned &|\Lambda^{-1/2}(v_n^\varepsilon(t\wedge \tau^n_{\varepsilon,M})-\theta_n^\varepsilon(t\wedge \tau^n_{\varepsilon,M} ))|^2+
2\varepsilon\kappa\int_0^{t\wedge\tau^n_{\varepsilon,M}} |\Lambda^{\alpha-\frac{1}{2}}(v_n^\varepsilon-\theta^\varepsilon_n)|^2ds
\\\leq &2\varepsilon\int_0^{t\wedge\tau^n_{\varepsilon,M}} \frac{\kappa}{2}|\Lambda^{\alpha-\frac{1}{2}}(v_n^\varepsilon-\theta^\varepsilon_n)|^2+C|\Lambda^{\alpha-\frac{1}{2}}v_n^\varepsilon|^2 ds
\\&+2\varepsilon\int_0^{t\wedge\tau^n_{\varepsilon,M}} \frac{\kappa}{2}|\Lambda^{\alpha-\frac{1}{2}}(\theta_n^\varepsilon-v_n^\varepsilon)|^2+ C\|v_n^\varepsilon\|_{L^p}^{N_0}|\Lambda^{-1/2}(\theta_n^\varepsilon-v_n^\varepsilon)|^2+C\|v^\varepsilon_n\|^4_{H^\delta}ds
\\&+2\sqrt{\varepsilon}\int_0^{t\wedge\tau^n_{\varepsilon,M}}\langle \Lambda^{-1/2}(v_n^\varepsilon-\theta^\varepsilon_n),\Lambda^{-1/2}(G(v_n^\varepsilon)-G(\theta^\varepsilon_n))dW\rangle
\\&+\varepsilon C\int_0^{t\wedge\tau^n_{\varepsilon,M}}|\Lambda^{-1/2}(v_n^\varepsilon-\theta^\varepsilon_n)|^2ds.\endaligned$$
Then Gronwall's lemma yields that
$$\aligned &|\Lambda^{-1/2}(v_n^\varepsilon(t\wedge \tau^n_{\varepsilon,M})-\theta_n^\varepsilon(t\wedge \tau^n_{\varepsilon,M} ))|^2
\leq \bigg[2C\varepsilon\int_0^{t\wedge\tau^n_{\varepsilon,M}} (|\Lambda^{\alpha-\frac{1}{2}}v_n^\varepsilon|^2+\|v^\varepsilon_n\|^4_{H^\delta}) ds\\
 &+ 2\sqrt{\varepsilon}|\int_0^{t\wedge\tau^n_{\varepsilon,M}}\langle \Lambda^{-1/2}(v_n^\varepsilon-\theta^\varepsilon_n),\Lambda^{-1/2}(G(v_n^\varepsilon)
 -G(\theta^\varepsilon_n))dW\rangle|\bigg]e^{\varepsilon C\int_0^{t\wedge\tau^n_{\varepsilon,M}}\|v_n^\varepsilon\|_{L^p}^{N_0}ds+Ct\varepsilon}
.\endaligned$$
Using (4.6) we have
$$\aligned &(E[\sup_{0\leq s\leq t\wedge\tau^n_{\varepsilon,M}} |\Lambda^{-1/2}(v_n^\varepsilon(s)-\theta_n^\varepsilon(s))|^2]^q)^{2/q}\\\leq&Ce^{\varepsilon CtM^{N_0/p}+Ct\varepsilon}\left[(\varepsilon Mt+\varepsilon M^2t)^2+q\varepsilon\int_0^t(E[\sup_{0\leq r\leq s\wedge\tau_{\varepsilon,M}^n} |\Lambda^{-1/2}(v_n^\varepsilon(r)-\theta_n^\varepsilon(r))|^2]^q)^{2/q}ds \right].\endaligned$$
By Gronwall's lemma we obtain that
$$\aligned &(E[\sup_{0\leq s\leq T\wedge\tau_{\varepsilon,M}^n} |\Lambda^{-1/2}(v_n^\varepsilon(s)-\theta_n^\varepsilon(s))|^2]^q)^{2/q}\\\leq&Ce^{\varepsilon CTM^{N_0/p}+CT\varepsilon}(\varepsilon MT+\varepsilon M^2T)^2\exp\left[ CqT\varepsilon e^{\varepsilon CTM^{N_0/p}+CT\varepsilon} \right].\endaligned$$
Fixing $M$ and taking $q=2/\varepsilon$ we have
$$\aligned &\varepsilon\log P(\sup_{0\leq t\leq T\wedge\tau_{\varepsilon,M}}\|\theta_n^\varepsilon(t)-v_n^\varepsilon(t)\|^2_{H^{-1/2}}>\eta)\\
\leq &\varepsilon\log \frac{E[\sup_{0\leq s\leq T\wedge\tau_{\varepsilon,M}} |\Lambda^{-1/2}(v_n^\varepsilon(s)-\theta_n^\varepsilon(s))|^{2q}]}{\eta^q}\\\leq & \log C(\varepsilon MT+\varepsilon M^2T)^2-2\log\eta+CTe^{\varepsilon CTM^{N_0/p}+CT\varepsilon}+\varepsilon CM^{N_0/p}+CT\varepsilon\\
\rightarrow& -\infty,\textrm{ as }\varepsilon\rightarrow0.\endaligned\eqno(4.20)$$
Therefore, there exists  $\varepsilon_0$ such that for every $\varepsilon$ satisfying $0<\varepsilon\leq\varepsilon_0$,
$$P(\sup_{0\leq t\leq T}\|\theta_n^\varepsilon(t)-v_n^\varepsilon(t)\|^2_{H^{-1/2}}>\eta,
\sup_{0\leq t\leq T}(\|v_n^\varepsilon(t)\|^2_{H^\delta}+\|v_n^\varepsilon(t)\|^p_{L^p})\leq M)\leq e^{-R/\varepsilon}.\eqno(4.21)$$
By Lemma 4.4 and (4.21), we know that there exists  $\varepsilon_0$ such that for every $\varepsilon$ satisfying $0<\varepsilon\leq\varepsilon_0$ we have
$$P(\sup_{0\leq t\leq T}\|\theta_n^\varepsilon(t)-v_n^\varepsilon(t)\|^2_{H^{-1/2}}>\eta)\leq2e^{-R/\varepsilon}.$$
Since $R$ is arbitrary, the desired result follows.
$\hfill\Box$

Now we can finish the proof of Theorem 4.2.
\vskip.10in

\noindent\textbf{Proof of Theorem 4.2}  By Lemmas 4.5 and 4.6, we have for every $R>0$ there exists  $N_2$ such that
$$P(\sup_{0\leq t\leq T}\|\theta_{N_2}^\varepsilon(t)-\theta^\varepsilon(t)\|^2_{H^{-1/2}}>\frac{\eta}{3})\leq e^{-R/\varepsilon} \ \textrm{ for any } \varepsilon\in (0,1];$$
and
$$P(\sup_{0\leq t\leq T}\|v_{N_2}^\varepsilon(t)-v^\varepsilon(t)\|^2_{H^{-1/2}}>\frac{\eta}{3})\leq e^{-R/\varepsilon}\ \textrm{ for any } \varepsilon\in (0,1].$$
 For such $N_2$, according to Lemma 4.7, there exists  $\varepsilon_0$ such that for every $\varepsilon$ satisfying $0<\varepsilon\leq\varepsilon_0$,
$$P(\sup_{0\leq t\leq T}\|\theta_{N_2}^\varepsilon(t)-v_{N_2}^\varepsilon(t)\|^2_{H^{-1/2}}>\frac{\eta}{3})\leq e^{-R/\varepsilon}.$$
Therefore, for every $\varepsilon$ satisfying $0<\varepsilon\leq\varepsilon_0$  we have
$$P(\sup_{0\leq t\leq T}\|\theta^\varepsilon(t)-v^\varepsilon(t)\|^2_{H^{-1/2}}>\eta)\leq 3e^{-R/\varepsilon}.$$
Since $R$ is arbitrary, we have
$$\lim_{\varepsilon\rightarrow0}\varepsilon\log P(\sup_{0\leq t\leq T}|\Lambda^{-1/2}(\theta^\varepsilon(t)-v^\varepsilon(t))|^2>\eta)=-\infty,$$
i.e. (4.10) holds. Hence  the proof of Theorem 4.2 is complete.  $\hfill\Box$

\section*{Appendix}

\th{Theorem A.1} Suppose that A.1)-A.3) hold. Then for any $\theta_0\in H^\delta\cap L^p$ with $p$ in Hypothesis 3.5 iii),  (3.8) has a unique solution
 $$\theta_v\in L^\infty([0,T],H^\delta\cap L^p)\cap L^2([0,T],H^{\delta+\alpha})\cap C([0,T],H^{-\beta})$$
and it has the following estimate:
$$\sup_{t\in[0,T]}(|\Lambda^\delta\theta_v(t)|^2+\|\theta_v(t)\|_{L^p}^p)+\int_0^T|\Lambda^{\delta+\alpha}\theta_v(s)|^2ds\leq C, \eqno(A.4)$$
where $C$ is some constant only depending on $|\Lambda^\delta\theta_0|,\|\theta_0\|_{L^p},T$ and $\int_0^T|v|_U^2ds$.

\proof In the following we will assume that $\delta<1$. The case for $\delta\geq1$ is similar.

[Step 1] We first establish the existence of solutions of the following
equation
$$\frac{d \theta(t)}{dt}+A_\alpha\theta(t)+w(t)\cdot \nabla\theta(t)= k_\sigma*G(\theta(t))v(t),\eqno(A.5)$$
 $$\theta(0)=\theta_0\in H^3$$
 with a given smooth function $w(t)$ which satisfies $div w(t)=0$ and $\sup_{t\in[0,T]}\|w(t)\|_{C^3}\leq C$.  Here $k_\sigma*G(\theta)$ means for $y\in U$, $k_\sigma*G(\theta)(y)=k_\sigma*(G(\theta)(y))$, where $k_\sigma$ is the periodic Poisson kernel in $\mathbb{T}^2$ given by $\widehat{k_\sigma}(\zeta)=e^{-\sigma|\zeta|},\zeta\in \mathbb{Z}^2$.

 Then we have the following apriori estimate
 $$\aligned \frac{d}{dt}|\Lambda^3 \theta|^2+2\kappa|\Lambda^{3+\alpha}\theta|^2\leq 2|\langle w\cdot\nabla\theta,\Lambda^{6}\theta\rangle|+2|\langle \Lambda^3\theta,\Lambda^3k_\delta*G(\theta)v\rangle|.\endaligned$$
 By Lemmas 2.1 and 2.2 we have that
$$|\langle \Lambda^{3-\alpha}(w\cdot\nabla \theta),\Lambda^{3+\alpha}\theta\rangle|\leq C\|w\|_{C^3(\mathbb{T}^2)}|\Lambda^{4-\alpha}\theta||\Lambda^{3+\alpha}\theta|\leq C|\Lambda^{3}\theta|^2+\kappa|\Lambda^{3+\alpha}\theta|^2,$$
where in the last inequality we use the interpolation inequality and Young's inequality.

Note that  we also have $$|\Lambda^3k_\sigma*G(\theta)v|\leq
C(\sigma)\|G(\theta)\|_{L_2(U,H)}|v|_U\leq C|v|_U(|\Lambda^\alpha
\theta|+1).$$ Thus,$$\aligned \frac{d}{dt}|\Lambda^3
\theta|^2+\kappa|\Lambda^{3+\alpha}\theta|^2\leq
C|v|_U(|\Lambda^3\theta|^2+1)+C|\Lambda^3\theta|^2.\endaligned$$ Then by the standard Galerkin
approximation we obtain that there exists a solution $\theta\in
L^\infty([0,T],H^3)\cap L^2([0,T],H^{3+\alpha})\cap C([0,T],H^1)$ of
(A.5).

\vskip.10in

[Step 2] Now we construct an approximation of (3.8).

We pick a smooth $\phi\geq0$, with supp $\phi\subset[1,2]$ and
$\int_0^\infty \phi=1$,  and for $\sigma>0$ let
$$U_\sigma[\theta](t):=\int_0^\infty \phi(\tau)(k_\sigma*R^\bot\theta)(t-\sigma\tau)d\tau,$$ where $k_\sigma$ is the periodic Poisson Kernel in
$\mathbb{T}^2$ given by $\widehat{k_\sigma}(\zeta)=e^{-\sigma|\zeta|},\zeta\in \mathbb{Z}^2$, and we set $\theta(t)=0$ for $t<0$.

We take a sequence $\delta_n\downarrow0$ and consider the equation
$$\frac{d \theta_n(t)}{dt}+A_\alpha\theta_n(t)+u_n(t)\cdot \nabla\theta_n(t)= k_{\delta_n}*G(\theta_n)v(t) \eqno(A.6)$$
with initial data $\theta_n(0)=k_{\delta_n}*\theta_0$ and
$u_n=U_{\delta_n}[\theta_n]$.

For a fixed $n$, this is a linear
equation in $\theta_n$ on each subinterval $[t^n_k,t^n_{k+1}]$ with
$t^n_k=k\delta_n$, since $u_n$ is determined by the values of
$\theta_n$ on the two previous subintervals.

By [Step 1], we obtain
the existence of a solution to (A.6) for fixed $n$. Moreover by
(A.1) the solution satisfies the following $L^p$ norm estimate:
$$\aligned \frac{d}{dt}\|\theta_n\|^p_{L^p}&=p\int|\theta_n|^{p-2}\theta_n(k_{\delta_n}*G(\theta_n)v-u_n\cdot\nabla \theta_n-\Lambda^{2\alpha}\theta_n)dx\\&\leq p\int|\theta_n|^{p-2}\theta_n k_{\delta_n}*G(\theta_n)vdx\leq C|v|_U(\|\theta_n\|_{L^p}^p+1).\endaligned$$
Here in the first inequality we use $div u_n=0$ and $\int|\theta_n|^{p-2}\theta_n\Lambda^{2\alpha}\theta_n\geq0$ (c.f. \cite[Lemma 3.2]{Re95}). Then Gronwall's lemma implies that$$\sup_{t\in[0,T]}\|\theta_n\|^p_{L^p}\leq (\|\theta_0\|_{L^p}^p+\int_0^T|v|_Udt)\exp{(\int_0^T|v|_Udt)}.$$
By (2.1) we have$$\sup_{t\in[0,T]}\|u_n\|^p_{L^p}\leq C(\|\theta_0\|_{L^p}^p+\int_0^T|v|_Udt)\exp{(\int_0^T|v|_Udt)}.$$
Now we prove the uniform $H^\delta$ estimate:
 $$\aligned \frac{d}{dt}|\Lambda^\delta \theta_n|^2+2\kappa|\Lambda^{\delta+\alpha}\theta_n|^2\leq 2|\langle \Lambda^{\delta}(u_n\cdot\nabla\theta_n),\Lambda^{\delta}\theta_n\rangle|+2|\langle \Lambda^\delta\theta_n,\Lambda^\delta k_{\delta_n}*G(\theta_n)v\rangle|.\endaligned$$
By \cite[Proposition 3.6]{Re95} we have
$$|\langle \Lambda^{\delta}(u_n\cdot\nabla\theta_n),\Lambda^{\delta}\theta_n\rangle|\leq \frac{\kappa}{4}|\Lambda^{\delta+\alpha}\theta_n|^2+\frac{\kappa}{8}|\Lambda^{\delta+\alpha}u_n|^2+C\|u_n\|_{L^p}^{N_0}|\Lambda^{\delta}\theta_n|^2+C\|\theta_n\|_{L^p}^{N_0}|\Lambda^{\delta}u_n|^2,$$
where $N_0=\frac{\alpha}{\alpha-\frac{1}{2}-\frac{1}{p}}$.

 By  A.2) we also obtain
$$|\langle \Lambda^\delta\theta_n,\Lambda^\delta k_{\delta_n}*G(\theta_n)v\rangle|\leq C|v|_U|\Lambda^{\delta}\theta_n|(|\Lambda^{\delta+\alpha}\theta_n|+1)\leq \varepsilon|\Lambda^{\delta+\alpha}\theta_n|^2+C(|v|^2_U|\Lambda^{\delta}\theta_n|^2+1).$$
Thus,
 $$\aligned |\Lambda^\delta \theta_n(t)|^2+\kappa\int_0^t|\Lambda^{\delta+\alpha}\theta_n|^2 ds \leq \int_0^t \left[ 2C\|u_n\|_{L^p}^{N_0}|\Lambda^{\delta}\theta_n|^2+C\|\theta_n\|_{L^p}^{N_0}|\Lambda^{\delta}u_n|^2
 +C(|v|^2_U|\Lambda^{\delta}\theta_n|^2+1) \right]ds.\endaligned$$
Note that
we have (here we cannot control $|\Lambda^\delta u_n|$ by $|\Lambda^\delta\theta_n|$ pointwisely in time)
$$\int_0^t|\Lambda^{\delta}u_n|^2 ds\leq C\int_0^t|\Lambda^{\delta}\theta_n|^2ds.$$
Using Gronwall's inequality and $L^p$ norm estimate above we obtain the uniform $H^\delta$ estimate for $\theta_n$.

Then by standard argument we know that $\theta_n$ converges to the solution $\theta_v$ of (3.8), which implies (A.4). The proof of uniqueness is the same as in \cite[Theorem 3.7]{Re95} by  A.3).
$\hfill\Box$

\vskip.10in

\th{Acknowledgement.} The authors would like to thank the referee for  many valuable  comments and suggestions. The authors also thank Rongchan Zhu for her helpful discussions.


\begin{thebibliography}{99}

\bibitem{BD98} M. Bou\'{e}, P. Dupuis,
A variational representation for certain functionals of Brownian
motion. \emph{Ann. Probab.}, \textbf{26}  (1998), No. 4, 1641-1659.



\bibitem{BD00} A. Budhiraja, P. Dupuis, A variational representation for positive functionals of infinite dimensional Brownian motion, \emph{Probab. Math. Statist.}, \textbf{20} (2000),39-61.

\bibitem{BDM08} A. Budhiraja, P. Dupuis, V. Maroulas, Large
 deviations for infinite
dimensional stochastic dynamical systems,  \emph{Ann. Probab.} \textbf{36} (2008), 1390-1420.

\bibitem{BY82} M. T. Barlow, M. Yor, Semi-martingale inequalities via the Garsia-Rodemich-Rumsey lemma, and applications to local time, \emph{J. Funct. Anal.}, \textbf{49} (1982),198-229.

\bibitem{CV06} L. Caffarelli, A. Vasseur, Drift diffusion equations with fractional diffusion and the quasi-geostrophic equation, \emph{Annals of Math.},
\textbf{171} (2010), No. 3, 1903-1930.

\bibitem{CM10}
I.~Chueshov, A.~Millet, {S}tochastic 2{D} hydrodynamical type systems:
  well posedness and large deviations, \emph{Appl. Math. Optim.} \textbf{61} (2010),
  379-420.

\bibitem{CMT94}P. Constantin, A. Majda, E. Tabak: Formation of strong fronts in the 2-D quasi-geostrophic thermal active scalar. \emph{Nonlinearity} \textbf{7} (1994), 1495-1533.

\bibitem{CW99} P. Constantin, J. Wu, Behavior of solutions of 2D quasi-geostrophic equations, \emph{SIAM J. Math. Anal.} \textbf{30} (1999), 937-948.

\bibitem{DZ92} G. Da Prato, J. Zabczyk, Stochastic Equations in Infinite Dimensions. Cambridge University Press, 1992.

\bibitem{Da76} B. Davis. On the $L^p$-norms of stochastic integrals and other martingales, \emph{Duke Math. J.}, \textbf{43} (1976),697-704.

\bibitem{DZ93} A. Dembo, O. Zeitouni. Large Deviations Techniques and Applications. Jones and Bartlett, Boston, 1993.

\bibitem{DM08} J. Duan,  A. Millet,
Large deviations for the Boussinesq equations under random influences,
\emph{Stochastic Process. Appl.} \textbf{119} (2009), 2052-2081.

\bibitem{Fl94} F. Flandoli, Dissipativity and invariant measures for stochastic Navier-Stokes equations, \emph{NoDEA}  \textbf{1} (1994), 403-423.

\bibitem{FG95} F. Flandoli, D. Gatarek, Martingale and stationary solutions for stochastic Navier-Stokes equations, \emph{Probability Theory and Related Fields} \textbf{102} (1995), 367-391.

\bibitem{FW84} M.I. Freidlin and A.D. Wentzell, \emph{Random
 perturbations of dynamical systems,} Translated from the Russian by
 Joseph Szu"cs.
 Grundlehren der Mathematischen Wissenschaften [Fundamental Principles
 of Mathematical Sciences], 260. Springer-Verlag, New York, 1984.

 \bibitem{GK96} I. Gy\"{o}ngy, N. Krylov, Existence of strong solutions for It\^{o}'s stochastic
equations via approximations, \emph{Probab. Theory Relat. Fields} \textbf{105} (1996), 143-158.


\bibitem{Ju04} N. Ju, Existence and Uniqueness of the Solution to the Dissipative 2D Quasi-Geostrophic Equations in the Sobolev Space, \emph{Communications in Mathematical Physics} \textbf{251} (2004), 365-376.

\bibitem{Ju05} N. Ju, On the two dimensional quasi-geostrophic equations, \emph{Indiana Univ. Math. J.} \textbf{54} No. 3 (2005), 897-926.

\bibitem{KNV07} A. Kiselev, F. Nazarov, A. Volberg,  Global well-posedness for the critical 2D dissipative quasi-geostrophic equation, \emph{Invent. math.} \textbf{167} (2007), 445-453.

\bibitem{Kr10} N.V. Krylov, It\^{o}'s formula for the $L_p$-norm of stochastic $W^1_p$-valued
processes, \emph{Probab. Theory Relat. Fields} \textbf{147} (2010), 583-605.


\bibitem{Li09} W. Liu, Large deviations for stochastic evolution equations with small
multiplicative noise, \emph{Appl. Math. Optim.} \textbf{61} (2010), 27-56.

  \bibitem{MSS09} U. Manna, S. S. Sritharan, P. Sundar,
Large deviations for the stochastic shell model of turbulence,
\emph{NoDEA Nonlinear Differential Equations Appl.} \textbf{16} (2009), 493-521.

\bibitem{On05} M. Ondrej\'{a}t, Brownian representations of cylindrical local martingales, martingale problem and strong markov property of weak solutions of spdes in Banach spaces, \emph{Czechoslovak Mathematical Journal} \textbf{55} (130)(2005), 1003-1039.

  \bibitem{P87}  J. Pedlosky: Geophysical Fluid Dynamics. New York: Springer-Verlag, 1987

\bibitem{PR07} C. Pr\'{e}v\^{o}t, M. R\"{o}ckner, A Concise Course on Stochastic Partial Differential Equations, Lecture Notes in Math., vol.1905, Springer, 2007.

 \bibitem{RZ08}  J. Ren, X. Zhang,  Freidlin-Wentzell's large deviations for stochastic evolution equations. \emph{J. Funct. Anal.} \textbf{254} (2008), 3148-3172.

\bibitem{Re95} S. Resnick, Danymical Problems in Non-linear Advective Partial Differential Equations, PhD thesis, University of Chicago, Chicago, 1995.

\bibitem{RZZ12} M. R\"{o}ckner, R.-C. Zhu, X.-C. Zhu,  Sub- and supercritical stochastic quasi-geostrophic equation, arXiv:1110.1984v4.

\bibitem{St70} E. Stein, Singular Integrals and Differentiability Properties of Functions, Princeton, NJ: Princeton
University Press, 1970.

    \bibitem{SS06} S.S. Sritharan, P. Sundar, Large deviations for the two-dimensional Navier-Stokes equations with multiplicative noise, \emph{Stoch. Proc. Appl.} \textbf{116} (2006), 1636-1659.

    \bibitem{St84} D.W.   Stroock, An Introduction to the Theory of Large Deviations, Springer, New York, 1984.

    \bibitem{SZ11} A. \'{S}wi\c{e}ch, J. Zabczyk, Large deviations for stochastic PDE with L\'{e}vy noise,
    \emph{J. Funct. Anal.} \textbf{260} (2011),  674-723.

\bibitem{Te84} R. Temam, Navier-Stokes Equations, North-Holland, Amsterdam, 1984.

\bibitem{Va66} S.R.S. Varadhan, Asymptotic probabilities and differential equations, \emph{Comm. Pure. Appl. Math.}. \textbf{19} (1966), 261-286.

\bibitem{Va67} S.R.S. Varadhan, Diffusion processes in small time intervals, \emph{Comm. Pure. Appl. Math.}. \textbf{20} (1967), 659-685.

\bibitem{XZ09} T.Xu, T.S. Zhang, On the small time asymptotics of the two-dimensional stochastic Navier-Stokes equations, \emph{Ann. Inst. H. Poincar\'{e} Probab. Statist.} \textbf{45} (4) (2009), 1002-1019.

\end{thebibliography}
\end{document}